\newif\ifsiamclass
\IfFileExists{siamart251216.cls}{%
  \siamclasstrue
  \documentclass[onefignum,onetabnum]{siamart251216}
}{%
  \siamclassfalse
  \documentclass[11pt]{article}
  \usepackage[margin=1in]{geometry}
}

\ifsiamclass\else
  \usepackage{amsthm}
\fi

\usepackage{amsmath,amsfonts,amssymb,mathtools}
\usepackage{bm}
\usepackage{booktabs}
\usepackage{graphicx}
\usepackage{xcolor}
\usepackage{enumitem}
\usepackage{algorithmic}
\usepackage{float}
\usepackage{booktabs,array,tabularx}
\usepackage{nicefrac}
\usepackage{xfrac}
\ifsiamclass\else
  \newenvironment{keywords}{\paragraph*{Keywords.}}{}
  \newenvironment{AMS}{\paragraph*{AMS subject classifications.}}{}
\fi
\providecommand{\headers}[2]{}

\makeatletter
\@ifundefined{theorem}{\newtheorem{theorem}{Theorem}[section]}{}
\@ifundefined{lemma}{}{}
\@ifundefined{proposition}{\newtheorem{proposition}[theorem]{Proposition}}{}
\@ifundefined{corollary}{}{}
\@ifundefined{remark}{}{}
\@ifundefined{definition}{}{}
\makeatother

\newcommand{\R}{\mathbb{R}}

\newcommand{\range}{\operatorname{range}}
\newcommand{\rank}{\operatorname{rank}}

\newcommand{\E}{\mathbb{E}}

\newcommand{\argmin}{\operatorname*{argmin}}

\newcommand{\norm}[1]{\left\|#1\right\|}
\newcommand{\ip}[2]{\left\langle #1,#2\right\rangle}

\definecolor{shikharnotecolor}{HTML}{7703FC}
\newif\ifshownotes
\shownotestrue

\ifshownotes
  \newcommand{\ssnote}[1]{\textcolor{shikharnotecolor}{\textbf{((}#1\textbf{))}}}
\else
  \newcommand{\ssnote}[1]{\ignorespaces}
\fi

\newcommand{\maybeincludegraphics}[2][]{%
  \IfFileExists{#2}{%
    \includegraphics[#1]{#2}%
  }{%
    \fbox{\parbox[c][0.28\textheight][c]{0.82\linewidth}{\centering Missing figure file:\\[0.5ex]\texttt{\detokenize{#2}}}}%
  }%
}
\headers{Noise-Adaptive Analog--Digital Preconditioning}
{S. Shah, R. Li, V. Kalantzis, Y. Xi}

\title{Noise-Adaptive Analog--Digital Preconditioning}

\ifpdf
\hypersetup{
	pdftitle={Noise-Adaptive Analog--Digital Preconditioning},
	pdfauthor={S. Shah, R. Li, V. Kalantzis, Y. Xi}
}
\fi

\author{
Shikhar Shah\thanks{Department of Mathematics, Emory University, Atlanta, GA, USA
(\email{shikhar.shah@emory.edu}, \email{yxi26@emory.edu}).}
\and
Rui Peng Li\thanks{Center for Applied Scientific Computing, Lawrence Livermore National
Laboratory, Livermore, CA, USA
(\email{li50@llnl.gov}).
{The work of S. Shah, R. Li, and Y. Xi was performed under the auspices of
the U.S. Department of Energy by Lawrence Livermore National Laboratory under
Contract DE-AC52-07NA27344 (LLNL-JRNL-2026513).}
}
\and
Vassilis Kalantzis\thanks{Thomas J. Watson Research Center, IBM, Yorktown Heights, NY, USA
(\email{vkal@ibm.com}).}
\and
Yuanzhe Xi\footnotemark[1]
}

\begin{document}

\maketitle

\begin{abstract}
Analog in-memory computing offers the potential to accelerate
approximate inverse preconditioning, but hardware errors can
degrade preconditioners designed in digital arithmetic. We develop
an analog--digital preconditioner for symmetric positive-definite
linear systems that combines analog block Jacobi with a small
digital low-rank residual correction. The correction space is
constructed by probing the same programmed analog arrays used in the
solve, targeting errors from both the coupling omitted by block
Jacobi and persistent hardware perturbations. Under a random
linear-operator model, we separate the programmed-mean error from
fresh read fluctuations, distinguish the effects of copy averaging
and repeated reads, and quantify the sampling accuracy of
columnwise probing. Under a centered additive write-noise model,
we also derive a computable regularization threshold and identify
the approximation error introduced by the regularization shift.
Experiments with an analog hardware simulator show that programming
errors can change the directions that block Jacobi reduces poorly.
In the tested comparisons with significant programming errors,
learning the correction space from the programmed copies reduces
FGMRES iteration counts relative to construction from digital
block Jacobi.
\end{abstract}



\section{Introduction}
\label{sec:introduction}

\begin{figure}[t]
    \centering
    \includegraphics[width=0.95\linewidth]{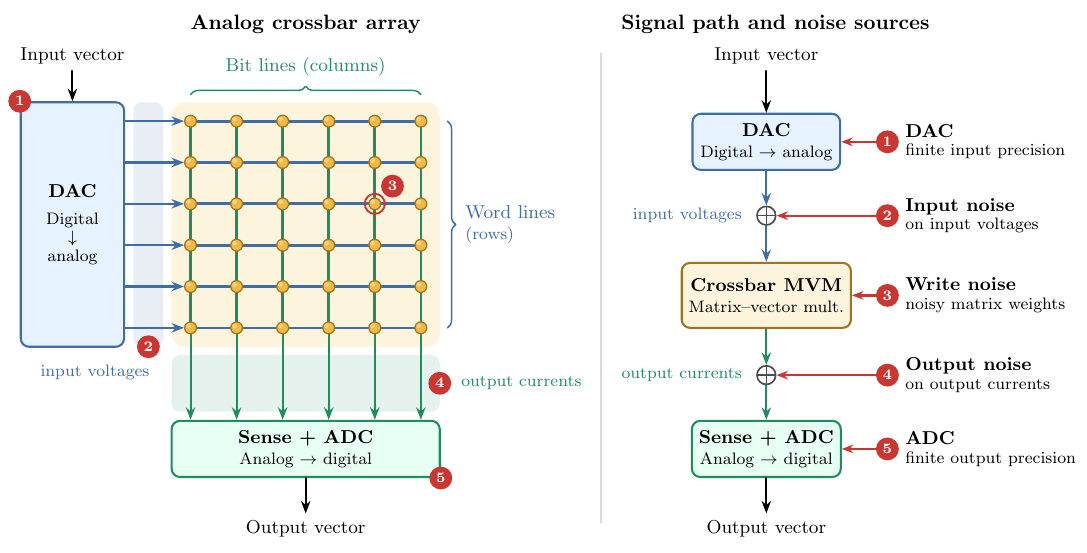}
    \caption{Analog matrix-vector multiplication on a crossbar array and the noise sources it introduces. Left: a DAC converts the digital input vector into voltages on the word lines, each cell stores one matrix entry as a conductance, and the currents collected on the bit lines are converted back to a digital output vector by the sense and ADC circuitry. Right: the same computation drawn as a signal path. Numbered markers correspond between the two panels and show where each error enters: (1) finite DAC precision, (2) noise on the input voltages, (3) write noise in the programmed conductances, (4) noise on the output currents, and (5) finite ADC precision.}
    \label{fig:crossbar_noise}
\end{figure}

Analog in-memory computing (AIMC) can perform matrix-vector multiplication efficiently by exploiting the physical properties of hardware, rather than executing individual floating-point
multiply-and-add instructions. For example, resistive crossbar arrays allow the implementation of matrix-vector multiplication by exploiting Ohm's law and Kirchhoff's laws, where matrix entries are encoded as conductances at crosspoint devices, while input vector elements are represented by applied voltages, as depicted in Figure~\ref{fig:crossbar_noise}. The resulting currents are accumulated along the crossbar columns, enabling many multiply-accumulate operations to be performed simultaneously \cite{hu2018memristor,li2018analogue}. The potential reduction in computational latency and energy consumption compared with conventional digital architectures
has motivated hybrid matrix algorithms that combine analog and digital computing
\cite{legallo2018mixed,sebastian2020memory}. On the other hand, analog matrix-vector multiplication is affected by programming inaccuracies,
device variability, and limited numerical
precision. Recent approaches aim to address these
challenges through error-aware programming,
calibration, and precision-extension techniques
\cite{buchel2023programming,kim2019precision,legallo2024fourquadrant}.\footnote{We note that, beyond electrical crossbars, photonic
architectures offer an alternative implementation
of matrix-vector products using optical
interference and programmable linear
transformations
\cite{shen2017deep,xie2025complex}.} Reviews of the devices, circuits, and architectures underlying AIMC are
given in \cite{Ielmini2018,sebastian2020memory,sensing}, and \cite{deeplearninganalog}
discusses analog accelerators from a systems perspective.

In this work, we focus on leveraging hybrid analog--digital preconditioning schemes
to solve the sparse linear system
\begin{equation}
    Ax=b,
\end{equation}
where $A \in \R^{N \times N}$ is symmetric positive-definite (SPD).
Since analog
hardware is an appealing architecture to perform matrix-vector multiplications, it has
long been recognized that analog computing can serve as a hardware accelerator of  approximate inverse preconditioners, e.g., \cite{benzi1999comparative,chen2021multiply,kalantzis2021solving,kalantzis2023fgmres}.
Analog in-memory hardware has been considered broadly for scientific computing since 2018. Feinberg et al.\ \cite{feinberg2018} showed how a fixed-point memristive accelerator can deliver the floating-point accuracy required by sparse iterative solvers, and Le Gallo et al.\ \cite{legallo2018mixed} combined an in-memory solver with high-precision digital refinement. Another demonstration is a memristor-based partial differential equation solver \cite{Zidan2018}. Hybrid analog--digital iterative linear system solvers typically perform matrix-vector
products with \(A\), inner products, residual checks, and vector addition/subtraction in digital arithmetic, while outsourcing the preconditioner application to analog arrays
\cite{feinberg2021analog,kalantzis2021solving,kalantzis2023fgmres,shah2026hybrid}.

Due to the nondeterministic noise introduced by each application of the approximate inverse preconditioner via analog hardware,
the outer loop is generally chosen to be a flexible iterative method, e.g.,
iterative refinement \cite{Kalantzis2025StableIS,legallo2018mixed,wu2026stable} or the flexible generalized minimal residual method (FGMRES) \cite{saad1993flexible}.
The effectiveness of this framework nevertheless depends
on the quality of the preconditioner delivered by the
programmed hardware \cite{feinberg2021analog}.
Block Jacobi preconditioning has shown comparatively robust
behavior under analog perturbations, with the additional benefit
that its dense block inverses can be applied independently
on separate arrays \cite{shah2026hybrid}. Its locality, however, leaves inter-block
and global coupling unresolved.

\subsection{Contributions}
In this paper, we aim to extend analog-based
block Jacobi preconditioning with a digital low-rank residual correction.
Low-dimensional corrections that complement local
preconditioners are well established in multigrid and
domain decomposition algorithms. For example, adaptive smoothed aggregation multigrid
constructs coarse information from the behavior of the
iteration itself \cite{brezina2005adaptive}, while
two-grid theory relates coarse-space quality to the
chosen relaxation
\cite{doi:10.1137/090772216, falgout2005twogrid,xu2017algebraic}.
GenEO solves local generalized eigenproblems to construct
robust coarse spaces \cite{spillane2014robust}, and
low-rank Schur-complement preconditioners capture coupling
omitted by local approximations
\cite{dillon2018hlr,li2016schur,xi2016multilevel}.
In the analog setting, the correction must account for
the programmed realization of the local preconditioner.
Programming errors can change which eigenmodes are
poorly treated, so a space constructed from digital
block Jacobi may miss important errors encountered
during the solve. Repeated reads reduce fresh fluctuations
around the programmed mean action but do not remove its
persistent discrepancy from the digital counterpart.
We therefore construct the correction space by probing
the same programmed copies used in the solve.
The setup combines forward analog applications with
digital matrix-vector products and orthogonalization;
the resulting basis remains fixed throughout the
Krylov iteration.

With a matrix \(Z\in\mathbb R^{N\times k}\)
of full column rank,
the corrected
preconditioner is
\begin{equation}
    \label{eq:intro_low_rank_correction}
    M_Z(B)
    =
    B+C_Z(I-AB),
    \qquad
    C_Z=Z(Z^TAZ)^{-1}Z^T.
\end{equation}
Here, \(B\) denotes one application of the programmed
analog block Jacobi preconditioner. For an input
residual \(r\), the analog action produces \(Br\),
and the digital correction \(C_Z(r-ABr)\) acts
on the remaining residual. The construction uses
the programmed action to identify the correction
space, allowing it to address both the local
approximation error and persistent hardware errors.

From \eqref{eq:intro_low_rank_correction}, $M_Z(B)$ can be rewritten as $(I - C_Z A)B + C_Z$, the adapted deflation two-level preconditioner with a smoother $B$. Deflation was introduced to accelerate conjugate gradients \cite{Nicolaides}, and two-level preconditioners derived from deflation, domain decomposition, and multigrid were later shown to be closely related \cite{Tang2009}. These works, however, assume a deterministic smoother and typically build $Z$ from $A$ or a digital preconditioner; in the present work, $B$ is a stochastic analog operator and $Z$ is learned from its programmed action.

Our analysis provides three results that guide
the construction and configuration of the preconditioner $M_Z(B)$:
\begin{enumerate}[label=(\roman*)]
    \item Under a random linear-operator model \eqref{eq:analog_action},
    we first derive a conditional mean--variance identity
    for one corrected application. It separates
    the squared programmed-mean error from fresh-read
    variance and shows how the digital correction
    affects both. Repeated reads reduce fresh variance,
    while an additional analysis under i.i.d.\ programming
    shows that copy averaging also affects programming
    variability.

    \item We analyze columnwise probing, with fresh
    reads for each vector application. Conditional
    mean and variance estimates show how each
    implemented sampling sweep approximates the
    programmed-mean error action. We also bound
    read-noise accumulation in an
    unorthogonalized sketch. These estimates
    quantify the effect of setup read averaging.

    \item Under a centered additive entrywise
    write-noise model, we derive an exact criterion
    for keeping the root-mean-square implementation
    error below the intended response for every
    fixed nonzero block input. The criterion gives
    a computable regularization threshold based on
    block spectral properties, hardware scaling,
    and copy count. We also identify the
    approximation error introduced by the shift.
\end{enumerate}

Simulator experiments compare correction spaces
constructed from digital block Jacobi with spaces
learned from the programmed copies, holding the
correction rank and programmed analog base fixed.
In the single-copy Poisson test, learning from the
programmed action reduces the mean FGMRES iteration
count from \(162.9\) to \(66.8\) at the default
noise level. Mean--variance diagnostics examine
which error contributions are reduced by correction
and averaging. Further experiments evaluate the
method on larger Poisson and interface Schur-complement
problems, while a finite element stiffness matrix
test compares the local regularization threshold
with the diagonal shifts giving the fewest FGMRES iterations.

Beyond the correction rank, the analog implementation offers two further
controls: averaging several independently programmed copies of the
preconditioner, and averaging repeated reads of each copy. All three improve
convergence at different hardware costs, so the preferred combination
depends on the array budget and energy model of a given architecture, which we
do not posit here.

The rest of the paper is organized as follows.
Section~\ref{sec:model} defines the analog block
Jacobi action and its digital low-rank correction.
Section~\ref{sec:theory} presents the mean--variance
analysis, probing construction, and regularization
criterion. Section~\ref{sec:experiments} reports
the numerical experiments, and
Section~\ref{sec:conclusion} concludes the paper.

\section{Analog Block Jacobi with Digital Low-Rank Correction}
\label{sec:model}

In this section we introduce analog block Jacobi
preconditioning with a digital low-rank correction.
Table~\ref{tab:notation} summarizes the
main notation for the preconditioner and its analysis.
Throughout the paper we use the energy inner product and its associated norm
\begin{equation}
    \ip{u}{v}_A=u^TAv,
    \qquad
    \norm{u}_A=\ip{u}{u}_A^{\nicefrac{1}{2}}.
\end{equation}

\begin{table}[tp]
    \centering
    \caption{Main notation for the analog--digital
    preconditioner and its analysis. The indices
    \(j,d,s\) identify a block, a programmed copy,
    and a read of that copy. In columnwise probing,
    \(t\) indexes the sampling sweep and \(a\) the
    probe column.}
    \label{tab:notation}
    \small
    \renewcommand{\arraystretch}{1.05}
    \setlength{\tabcolsep}{5pt}
    \begin{tabularx}{\linewidth}{
        @{}l >{\raggedright\arraybackslash}X@{}
    }
        \toprule
        \textbf{Symbol} & \textbf{Meaning} \\
        \midrule
        \multicolumn{2}{@{}l}{
            \textit{Digital target and implementation parameters}
        }\\
        \addlinespace[0.3ex]

        \(p\)
        & Number of blocks and arrays per programmed copy. \\

        \(D\)
        & Number of independently programmed copies. \\

        \(S\)
        & Repeated reads per copy within one averaged
          vector application. \\

        \(\beta\)
        & Nonnegative regularization shift. \\

        \(R_j,\ A_j\)
        & Restriction to block \(j\) and local matrix
          \(A_j=R_jAR_j^T\). \\

        \(M_{\beta,p},\ M_p\)
        & Regularized digital block Jacobi target and
          its unregularized counterpart,
          \(M_p=M_{0,p}\). \\

        \midrule
        \multicolumn{2}{@{}l}{
            \textit{Programming and read randomness}
        }\\
        \addlinespace[0.3ex]

        \(W_d,\ W\)
        & Programming realization of copy \(d\) and
          collection of all \(D\) realizations;
          fixed during setup and solve. \\

        \(\omega_{ds},\ \omega_d\)
        & Randomness of read \(s\) of copy \(d\),
          and of a generic read of that copy. \\

        \(\omega\)
        & Collection of fresh randomness in
          \(DS\) reads of one vector application. \\

        \(\omega_{t,a}\)
        & Fresh read collection for probe column
          \(a\) at sweep \(t\), comprising \(S\) reads
          from each of the same \(D\) programmed copies. \\

        \(B(W_d,\omega_{ds})\)
        & Analog operator for read $s$ of copy \(d\). \\

        \(B(W,\omega)\)
        & Analog operator averaged over \(D\) copies
          and \(S\) reads per copy, all applied to
          the same input vector. \\

        \(\overline B_{W_d}\)
        & Conditional mean action of copy \(d\),
          with \(W_d\) fixed. \\

        \(\overline B_W\)
        & Conditional mean of the averaged action,
          with \(W\) fixed. \\

        \(\overline E_W\)
        & Error operator \(I-\overline B_WA\) of programmed mean
          action before digital correction. \\

        \(\Xi_d\)
        & Persistent discrepancy
          \(\overline B_{W_d}-M_{\beta,p}\). \\

        \(\overline\Xi_D\)
        & Average persistent discrepancy,
          \(\frac{1}{D}{\sum_{d=1}^D\Xi_d}\). \\

        \(\Delta_d(\omega_d)\)
        & Fresh fluctuation of one read about
          \(\overline B_{W_d}\). \\

        \(\Delta_W(\omega)\)
        & Averaged fresh fluctuation,
          \(B(W,\omega)-\overline B_W\). \\

        \midrule
        \multicolumn{2}{@{}l}{
            \textit{Digital correction}
        }\\
        \addlinespace[0.3ex]

        \(r,\ z,\ v\)
        & Residual, computed and exact corrections:
          \(z=M_Z(B)r\),
          \(v=A^{-1}r\). \\

        \(Z\)
        & Full-column-rank $N \times k$ correction basis matrix, fixed
          during the solve.\\

        $k$
        &
        Rank of $Z$.
    If selecting \(k\) would split a complex-conjugate Schur pair, both vectors are retained, increasing the actual correction rank to \(k+1\).\\

        $q$    & Number of sampling sweeps that build the correction basis. \\
        $\ell$ & Probe width: number of probe columns in each sampling sweep. \\
    
        \(A_Z\)
        & Reduced SPD matrix \(Z^TAZ\). \\

        \(C_Z\)
        & Residual-to-correction map
          \(ZA_Z^{-1}Z^T\). \\

        \(P_Z\)
        & \(A\)-orthogonal projector \(C_ZA\)
          onto \(\range(Z)\). \\

        \(M_Z(B)\)
        & Complete analog--digital preconditioner,
          \(B+C_Z(I-AB)\). \\

        \bottomrule
    \end{tabularx}
\end{table}

\subsection{Multi-copy analog block Jacobi base}

Let \(\{\mathcal I_j\}_{j=1}^p\) be a nonoverlapping
partition of the unknowns, and
\(R_j\in\R^{|\mathcal I_j|\times N}\)
be the restriction operator
to block \(\mathcal I_j\). The exact-arithmetic
block Jacobi approximate inverse is
\begin{equation}
    \label{eq:block_jacobi}
    M_p=\sum_{j=1}^p R_j^T A_j^{-1}R_j,
    \qquad
    A_j=R_jAR_j^T.
\end{equation}
For a regularization shift \(\beta\geq0\), define
\begin{equation}
    \label{eq:regularized_block_jacobi}
    M_{\beta,p}
    =
    \sum_{j=1}^p R_j^T(A_j+\beta I)^{-1}R_j,
\end{equation}
so that \(M_{0,p}=M_p\).
The shift controls the sensitivity of the local
inverse actions to programming errors.
Section~\ref{sec:regularization} gives a
signal-to-noise criterion for selecting \(\beta\)
under an additive write-noise model.
Suppose that each local inverse \((A_j+\beta I)^{-1}\) is
implemented on a separate analog array and is independently applied
to the corresponding subvector \(R_jr\) of an
input \(r\in\R^N\).
The \(p\) analog arrays together
form one programmed copy of the target
\(M_{\beta,p}\). We use \(D\) independently
programmed copies, requiring \(pD\) arrays.

Let \(W_{d,j}\) denote the programming realization
of the array storing block \(j\) in copy~\(d\).
It includes write errors and other hardware
distortions that remain \emph{fixed} after programming.
We collect these realizations as
\begin{equation}
    W_d=(W_{d,1},\ldots,W_{d,p}),
    \qquad
    W=(W_1,\ldots,W_D),
\end{equation}
where \(W_d\), \(d=1, 2, \ldots, D\), describes one complete copy, and
\(W\) describes all \(D\) copies.
One read of a copy applies its \(p\) arrays to
the corresponding input subvectors and assembles
their outputs into a vector in \(\R^N\).
Let \(S\) be the number of reads of each
copy within a preconditioner application,
all using the same input \(r\).
For 
\(s=1,\ldots,S\),
let \(\omega_{ds}\) denote the random variable modeling the read-dependent effects
of read \(s\) of copy \(d\).
For the analysis, we represent this read by a
random linear operator:
\begin{equation}
    \label{eq:analog_action}
    B(W_d,\omega_{ds})r
    =
    M_{\beta,p}r+\delta_{W_d,\omega_{ds}}(r),
\end{equation}
where
the perturbation \(\delta_{W_d,\omega_{ds}}(r)\)
collects both persistent programming discrepancies
and read-dependent fluctuations across the
\(p\) arrays, and is linear in \(r\) for fixed
$(W_d, \, \omega_{ds})$. Since each array acts only
on its corresponding block,
\(B(W_d,\omega_{ds})\) is also block diagonal with
the same partition as \(M_{\beta,p}\).
This is an analysis model: quantization, clipping,
input-dependent scaling, and input-independent
additive read noise cannot, in general, be represented
exactly by a random linear operator.

\begin{figure}[t]
    \centering
    \includegraphics[width=0.95\linewidth]{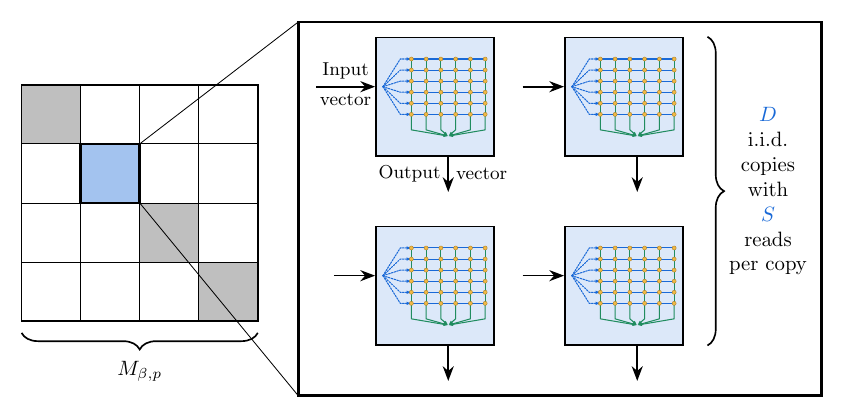}
    \caption{Multi-copy implementation of
    \eqref{eq:averaged_action}.
    Each copy uses $p$ arrays, one per block.
    The \(D\) copies are programmed independently,
    each is read \(S\) times with the same input,
    and their outputs are averaged digitally.
    The inset shows the \(D\) arrays storing
    block \(2\). If all \(pD\) arrays can be
    read in parallel, one application requires
    \(S\) sequential read rounds and \(pDS\)
    individual array reads. With constant
    per-array read latency and energy, these
    counts give \(\mathcal O(S)\) analog-read
    time and \(\mathcal O(pDS)\) read energy.
    These estimates exclude digital averaging,
    data movement, and the low-rank correction.}
    \label{fig:double_avg}
\end{figure}

Let \(\omega=\{\omega_{ds}\}_{d=1,\ldots,D;\,s=1,\ldots,S}\)
collect the read randomness in one averaged application.
The analog preconditioner applies all reads to the same
input \(r\) and averages their outputs digitally:
\begin{equation}
    \label{eq:averaged_action}
    B(W,\omega)r
    :=
    \frac1{DS}\sum_{d=1}^D\sum_{s=1}^S
        B(W_d,\omega_{ds})r.
\end{equation}
The counts $D$ and $S$ are specified for each configuration
and suppressed in the notation \(B(W,\omega)\).
The programming realizations \(W\) remain fixed
throughout setup and solve, whereas \(\omega\)
is resampled for each application. Figure~\ref{fig:double_avg}
illustrates the implementation and its idealized
read costs.

\subsection{Programmed mean and read fluctuations}
\label{sec:programmed_mean}

Repeated reads of a fixed programmed copy fluctuate
around a mean action that may differ from the digital
target \(M_{\beta,p}\). Within the random linear-operator
model \eqref{eq:analog_action},
we separate this persistent discrepancy from
read fluctuations and describe how both enter
the averaged preconditioner.
Let \(\omega_d\) denote the read randomness in a
generic read of copy \(d\). To ensure that the mean
actions and mean-square fluctuations are
finite, we assume that
\begin{equation}
    \E\!\left[
        \norm{B(W_d,\omega_d)}_F^2
    \right]<\infty,
    \qquad d=1,\ldots,D,
\end{equation}
where \(\norm{\cdot}_F\) is the Frobenius norm
and the expectation is over programming and read
randomness. For a fixed programming realization
\(W_d\), define the mean action
\begin{equation}
    \label{eq:device_mean}
    \overline B_{W_d}
    :=
    \E_{\omega_d}[B(W_d,\omega_d)\mid W_d].
\end{equation}
This conditional expectation averages over
reads while holding the programmed copy fixed.
Define its discrepancy from the digital target by
\begin{equation}
\label{eq:Xi_d}
    \Xi_d:=\overline B_{W_d}-M_{\beta,p},
\end{equation}
and the fluctuation of an individual read by
\begin{equation}
    \label{eq:device_fluctuation}
    \Delta_d(\omega_d)
    :=
    B(W_d,\omega_d)-\overline B_{W_d}.
\end{equation}
Because \(\overline B_{W_d}\) is the conditional
mean of the read operator, this fluctuation has
zero conditional mean:
\begin{equation}
    \label{eq:device_fluctuation_mean}
    \E_{\omega_d}[\Delta_d(\omega_d)\mid W_d]=0.
\end{equation}
One read therefore decomposes as
\begin{equation}
    \label{eq:three_layer}
    B(W_d,\omega_d)
    =
    M_{\beta,p}+\Xi_d+\Delta_d(\omega_d).
\end{equation}
For an input \(r\),
we have
\(\delta_{W_d,\omega_{ds}}(r)=\Xi_dr+\Delta_d(\omega_{ds})r\)
for the perturbation
in \eqref{eq:analog_action}.
The discrepancy \(\Xi_d\) remains fixed across
reads.

We now pass from a single copy to the full
collection \(W=(W_1,\ldots,W_D)\) and the averaged
operator \(B(W,\omega)\).
For the averaging analysis, we assume that,
conditional on \(W\), read randomness is
independent across copies, reads, and
preconditioner applications.
For each copy \(d\), the conditional distribution
of \(\omega_{ds}\) given \(W\) is the same as that
of \(\omega_d\) given \(W_d\).
The single-read
distributions and read settings remain unchanged
when \(S\) varies. These assumptions describe
a fixed operating regime without drift or
conditional correlations between distinct reads.
We call read randomness with these properties \textit{fresh}: unlike the programming realization \(W\), it is drawn independently for each read.
The programming realizations \(W_1,\ldots,W_D\)
are independent, as specified by the multi-copy
construction. Dependence among programming errors
within one copy, or among the \(p\) array
fluctuations within a single read, is allowed.
Section~\ref{sec:regularization} introduces
additional assumptions for its write-noise model.

The persistent discrepancy and fresh fluctuation
of the averaged action are
\begin{equation}
    \overline\Xi_D
    :=
    \frac1D\sum_{d=1}^D\Xi_d,
    \qquad
    \Delta_W(\omega)
    :=
    \frac1{DS}\sum_{d=1}^D\sum_{s=1}^S
        \Delta_d(\omega_{ds}).
\end{equation}
Substituting the single-copy decomposition
\eqref{eq:three_layer} into
\eqref{eq:averaged_action} gives
\begin{equation}
    \label{eq:Bavg}
    B(W,\omega)
    =
    M_{\beta,p}+\overline\Xi_D+\Delta_W(\omega).
\end{equation}
Each \(\Xi_d\) remains unchanged across the
\(S\) reads of copy \(d\), so the persistent
discrepancy is averaged only across copies.
The fresh fluctuations are averaged over all
\(DS\) reads.

Under the conditional read-distribution assumption,
\eqref{eq:device_fluctuation_mean} gives
\(\E_\omega[\Delta_d(\omega_{ds})\mid W]=0\)
for every \(d,s\). Therefore, their average
\(\Delta_W(\omega)\) also has zero
conditional mean. Taking the conditional
expectation in \eqref{eq:Bavg}, with
\(\overline\Xi_D\) fixed by \(W\), yields
\begin{equation}
    \label{eq:programmed_mean_smoother}
    \begin{aligned}
        \overline B_W
        &:=\E_\omega[B(W,\omega)\mid W]\\
        &=M_{\beta,p}+\overline\Xi_D
          =\frac1D\sum_{d=1}^D\overline B_{W_d}.
    \end{aligned}
\end{equation}
Thus, the averaged action separates into its
programmed mean and fresh fluctuation:
\begin{equation}
    \label{eq:conditional_decomposition}
    B(W,\omega)=\overline B_W+\Delta_W(\omega),
    \qquad
    \E_\omega[\Delta_W(\omega)\mid W]=0.
\end{equation}
In particular,
\(\Delta_W(\omega)=B(W,\omega)-\overline B_W\).

\subsection{Digital low-rank residual correction}
For a fixed input residual \(r\), the preconditioner
returns an approximate correction \(z\approx A^{-1}r\).
The digital correction lies in \(\range(Z)\), where
\(Z\in\R^{N\times k}\) has full column rank and
\(k\) denotes the actual correction rank.
The basis is constructed by probing the same programmed
copies used in the solve
(Section~\ref{sec:learn_programmed_space}) and remains
fixed throughout the Krylov iteration.

For one averaged analog application \(B=B(W,\omega)\),
the initial approximation and its remaining residual are
\begin{align}
    z_0=Br, \qquad
    s_0 = r-Az_0=(I-AB)r.
    \label{eq:unresolved_residual}
\end{align}
We choose \(z=z_0+Zc\) to minimize the energy-norm
error:
\begin{equation}
    c=\argmin_{\xi\in\R^k}
        \norm{z_0+Z\xi-A^{-1}r}_A^2.
\end{equation}
The first-order optimality condition is
\(Z^TA(z-A^{-1}r)=0\), or equivalently
\(Z^T(r-Az)=0\). Thus, the corrected residual
is orthogonal to \(\range(Z)\).
Substituting \(r-Az=s_0-AZc\) gives the reduced
system \(Z^TAZc=Z^Ts_0\). Define
\begin{equation}
    \label{eq:low_rank_correction_operator}
    A_Z=Z^TAZ,
    \qquad
    C_Z=ZA_Z^{-1}Z^T.
\end{equation}
Since \(A\) is SPD and \(Z\) has full column rank,
\(A_Z\) is SPD. Solving for \(c\) gives
\begin{equation}
    \label{eq:low_rank_corrected_action}
    z=z_0+Zc=z_0+C_Zs_0=M_Z(B)r,
\end{equation}
where
$M_Z(B):=B+C_Z(I-AB)$.
All operations after \(Br\) are digital; neither
\(C_Z\) nor \(M_Z(B)\) is explicitly assembled.

To describe the effect on the approximation error,
let \(v=A^{-1}r\) denote the exact correction and
define \(P_Z=C_ZA=Z(Z^TAZ)^{-1}Z^TA\), the
\(A\)-orthogonal projector onto \(\range(Z)\).
Since \(s_0=r-Az_0=A(v-z_0)\), we have
\begin{equation}
        C_Zs_0
        =C_ZA(v-z_0)=P_Z(v-z_0),
\end{equation}
from which it follows that
\begin{equation}
            v-z
        =(v-z_0)-C_Zs_0
          =(I-P_Z)(v-z_0).
\end{equation}
The digital correction removes the
\(A\)-orthogonal projection of the error left
by the analog approximation.
Substituting \(z_0=BAv\) and \(z=M_Z(B)Av\)
gives \((I-M_Z(B)A)v=(I-P_Z)(I-BA)v\).
Since this holds for every \(v\), each realized
\(B\) satisfies
\begin{equation}
    \label{eq:corrected_error_identity}
    I-M_Z(B)A=(I-P_Z)(I-BA).
\end{equation}
Section~\ref{sec:theory} uses this identity to
analyze the mean-square error of one corrected
preconditioner application.

\section{Analysis and Design of the Analog--Digital Preconditioner}
\label{sec:theory}

In this section, we perform an analysis to identify which error sources are addressed
by digital correction, read averaging, copy averaging,
and local regularization. The stated identities and bounds
hold under their respective models, and they do not jointly
provide a convergence guarantee or an optimal parameter
choice for the complete implementation. We use them to
motivate the construction and identify candidate
configurations, whose solver performance is experimentally assessed
in Section~\ref{sec:experiments}.

\subsection{Programmed mean error and fresh variance}
When an analog preconditioner produces an inaccurate
correction, repeated reads, additional programmed copies,
and a different correction space address different
sources of that error. We first analyze one application to
distinguish these effects and identify which errors
read averaging can reduce.

For a fixed input \(r\), the desired correction is
\(v=A^{-1}r\), whereas the implemented preconditioner
returns \(M_Z(B)r=M_Z(B)Av\).
Their difference, \(v-M_Z(B)Av\), measures the
combined approximation and hardware errors remaining
after digital correction. 
During the solve, the programmed copies \(W\) and
correction basis \(Z\) remain fixed, while read
fluctuations vary between applications. We therefore
hold \(W,Z\) fixed and examine the expected squared
\(A\)-norm of the application error over fresh reads.

Define the programmed-mean error operator by
\begin{align}
    \overline E_W
    =I-\overline B_WA =(I-M_{\beta,p}A)-\overline\Xi_DA.
    \label{eq:mean_error_split}
\end{align}
It combines the approximation error of digital block
Jacobi with the persistent programming discrepancy.
Substituting \(B=B(W,\omega)=\overline B_W+\Delta_W(\omega)\)
into \eqref{eq:corrected_error_identity} gives
\begin{align}
    v-BAv
    &=\overline E_Wv-\Delta_W(\omega)Av,
    \label{eq:realized_error_decomposition}\\
    v-M_Z(B)Av
    &=(I-P_Z)\overline E_Wv
      -(I-P_Z)\Delta_W(\omega)Av.
    \label{eq:two_grid_split_clean}
\end{align}

\begin{theorem}[Conditional mean--variance identity]
\label{thm:one_step_bias_variance_clean}
Assume the read model in
Section~\ref{sec:programmed_mean}.
Conditional on \(W\), assume that the fresh reads
in the current application are independent of all
randomness used to construct \(Z\).
For every fixed \(v\in\R^N\), write \(B=B(W,\omega)\)
and define
\begin{equation}
    \begin{aligned}
        e(\omega)&:=\left(I-M_Z(B)A\right)v,\\
        \bar e&:=(I-P_Z)\overline E_Wv,\\
        \eta_d(\omega_d)&:=(I-P_Z)\Delta_d(\omega_d)Av.
    \end{aligned}
\end{equation}
These vectors denote the corrected application error,
its programmed-mean component, and the projected
single-read fluctuation, respectively.
Then
\begin{align}
    \E_\omega[e(\omega)\mid W,Z]&=\bar e,
    \label{eq:mean_corrected_error_identity}\\
    \E_\omega[\norm{e(\omega)}_A^2\mid W,Z]
    &=\norm{\bar e}_A^2
    +\frac1{D^2S}\sum_{d=1}^D
    \E_{\omega_d}[\norm{\eta_d(\omega_d)}_A^2\mid W,Z].
    \label{eq:one_step_bias_variance_clean}
\end{align}
In \eqref{eq:one_step_bias_variance_clean}, the
left-hand side is the conditional mean-square
application error. The first term on the right-hand side is
the squared conditional mean error, and the second
is the conditional read variance, both measured
in the \(A\)-norm.
The expectations are over fresh reads with \(W,Z\)
fixed; \(\omega_d\) denotes a generic fresh read
of copy \(d\).
\end{theorem}

\begin{proof}
Define the averaged projected fluctuation by
\begin{equation}
    \eta(\omega):=\frac1{DS}
        \sum_{d=1}^D\sum_{s=1}^S\eta_d(\omega_{ds}).
\end{equation}
Equation~\eqref{eq:two_grid_split_clean} becomes
\(e(\omega)=\bar e-\eta(\omega)\).
The independence assumption ensures that conditioning
on \(Z\) leaves the current read distribution unchanged.
The single-read fluctuations therefore remain centered
and independent given \(W,Z\), so
\(\E_\omega[\eta(\omega)\mid W,Z]=0\).
Taking conditional expectations gives
\eqref{eq:mean_corrected_error_identity}.
The cross term between \(\bar e\) and \(\eta(\omega)\)
also vanishes, yielding
\begin{equation}
    \label{eq:expanded_square}
    \E_\omega[\norm{e(\omega)}_A^2\mid W,Z]
    =\norm{\bar e}_A^2
    +\E_\omega[\norm{\eta(\omega)}_A^2\mid W,Z].
\end{equation}
Expanding the squared norm of the average
\(\eta(\omega)\), conditional independence and
centering eliminate cross terms between distinct reads.
For each copy, the \(S\) remaining squared terms
have identical conditional expectations, giving
\begin{equation}
    \E_\omega[\norm{\eta(\omega)}_A^2\mid W,Z]
    =\frac1{D^2S}\sum_{d=1}^D
        \E_{\omega_d}[\norm{\eta_d(\omega_d)}_A^2\mid W,Z].
\end{equation}
Here, the coefficient is \(S/(DS)^2=1/(D^2S)\).
Substitution into \eqref{eq:expanded_square} proves
\eqref{eq:one_step_bias_variance_clean}.
\end{proof}
Equation~\eqref{eq:one_step_bias_variance_clean}
shows how additional reads of the existing programmed
copies affect the mean-square correction error.
For fixed \(W,Z,v\), and single-read
distributions, increasing \(S\) reduces the second
term on the right-hand side
in \eqref{eq:one_step_bias_variance_clean}
as \(1/S\), but leaves the first term unchanged.
Thus, additional reads are useful for suppressing
fresh fluctuations.

The read-variance term in
\eqref{eq:one_step_bias_variance_clean} equals the
average projected single-read variance across copies,
divided by \(DS\). Using \(\eta_d(\omega_d)\) from
Theorem~\ref{thm:one_step_bias_variance_clean}, this term is
\begin{equation}
\label{eq:read-var}
    \frac1{DS}
    \left[
        \frac1D\sum_{d=1}^D
        \E_{\omega_d}[\norm{\eta_d(\omega_d)}_A^2\mid W,Z]
    \right].
\end{equation}
For example, if each copy has the same projected
single-read variance \(\sigma^2\), this term
is \(\sigma^2/(DS)\).
For configurations with the same average variance, \eqref{eq:read-var}
scales exactly as
\(1/(DS)\). Adding copies or relearning \(Z\) may
change the average.

Copy averaging also affects the persistent discrepancy
\(\overline\Xi_D\) in the programmed mean, as shown
by \eqref{eq:programmed_mean_smoother}, while
repeated
reads do not.
To quantify this effect, we now average over
programming realizations, keeping the digital target
\(M_{\beta,p}\) and the input direction \(v\) fixed.
Assume additionally that \(\Xi_1,\ldots,\Xi_D\)
are independent and identically distributed with
finite second moments.
Write \(\E_W\) for expectation over the programming
realizations \(W\), and let \(\mu_\Xi := \E_W[\Xi_1]\) denote the common mean of the \(\Xi_d\). For quantities involving the
learned basis \(Z\), write \(\E\) for the joint
expectation over programming realizations and the
random probes and setup reads used to construct \(Z\).

We now examine the first term of
\eqref{eq:one_step_bias_variance_clean}.
The \(A\)-orthogonal projection yields
\begin{equation}
\label{eq:first_term_A_norm}
    \E\norm{(I-P_Z)\overline E_Wv}_A^2
    =
    \E_W\norm{\overline E_Wv}_A^2
    -
    \E\norm{P_Z\overline E_Wv}_A^2,
\end{equation}
where we first evaluate the squared programmed-mean
error before correction and then subtract its component
projected by \(Z\).

Then, we have
\begin{equation}
    \label{eq:mean_error_copy_averaging}
    \begin{aligned}
        \E\norm{(I-P_Z)\overline E_Wv}_A^2
        =&\norm{(I-M_{\beta,p}A-\mu_\Xi A)v}_A^2
        +\frac1D\E_W\norm{(\Xi_1-\mu_\Xi)Av}_A^2\\
        &-\E\norm{P_Z\overline E_Wv}_A^2.
    \end{aligned}
\end{equation}

\begin{figure}[t]
    \centering
    \includegraphics[width=0.95\linewidth]{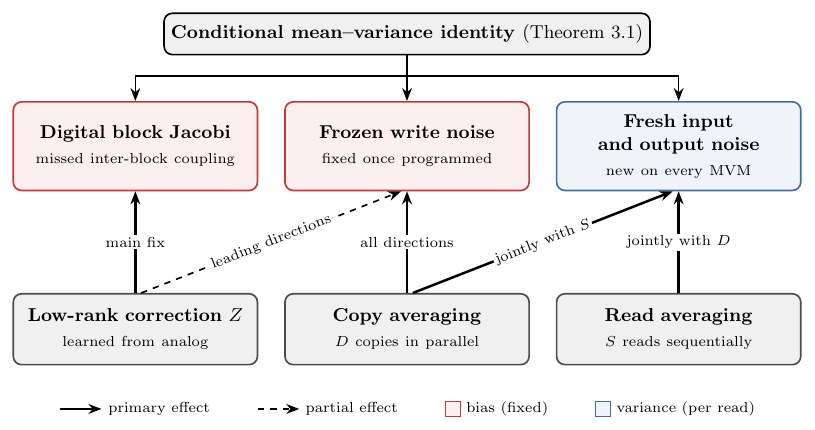}
    \caption{Effects of digital correction and averaging.
    Digital correction removes the \(A\)-orthogonal
    projection of both the programmed-mean error and
    fresh read fluctuations onto \(\range(Z)\),
    leaving the components projected by \(I-P_Z\).
    For fixed input, programmed copies, basis, and
    single-read distributions, increasing \(S\)
    reduces fresh-read variance as \(1/S\), without
    changing the programmed-mean error.
    Under independent and identically distributed programming,
    with the digital target and single-copy distribution
    fixed, averaging \(D\) copies reduces mean-square
    programming variability before correction as \(1/D\),
    but retains any common programming bias.}
    \label{fig:error_decomp}
\end{figure}

For a fixed digital target, input direction, and
single-copy distribution, the first term on the
right-hand side is unchanged by copy averaging: it combines
the digital approximation error with the common
programming bias. The second term decreases as
\(1/D\), showing that additional copies reduce
programming variability that repeated reads
cannot remove. The final term accounts for
the mean-error component removed by \(Z\).
Because this captured component also depends on
the programmed copies and learned basis, the
\(1/D\) reduction alone does not guarantee
monotonic improvement of the corrected mean error.
Figure~\ref{fig:error_decomp} summarizes these
distinct effects of averaging and correction.

\subsection{Learning the correction space from the programmed realization}
\label{sec:learn_programmed_space}

We construct the correction space from the errors observed
on the same programmed copies used in the solve.
%
Theorem~\ref{thm:one_step_bias_variance_clean}
shows that digital correction leaves the mean-error
component \((I-P_Z)\overline E_Wv\).
For a prescribed rank of the subspace correction, we aim to reduce
the worst-case remaining mean error, measured by the
operator $A$-norm 
\begin{equation}
    \label{eq:oracle_correction_space}
    \norm{(I-P_Z)\overline E_W}_A
    =
    \sup_{\norm{v}_A=1}
    \norm{(I-P_Z)\overline E_Wv}_A.
\end{equation}
The correction space \(\range(Z)\) should therefore
capture the largest mean-error outputs \(\overline E_Wv\),
measured in the \(A\)-norm. These outputs include both
digital approximation errors and persistent programming
discrepancies. We construct candidate correction
directions by randomized probing.
We use a forward-only probing scheme inspired by
randomized range finding \cite{halko2011finding}.
Let \(k\) be the requested correction rank,
\(k\leq\ell\leq N\) the probe width, and \(q\geq1\)
the number of sampling sweeps.
Choose a random starting block
\(\Omega\in\R^{N\times\ell}\) and set \(Y_0=\operatorname{orth}(\Omega)\).
Conditional on \(W\), this block is independent
of fresh setup reads.

We first build a candidate
correction space for the
programmed-mean error operator \(\overline E_W\)
by subspace iteration.
Each sweep samples its action on the current basis vectors
and orthonormalizes the resulting vectors.
At sweep \(t=1,\ldots, q\), holding
\(Y_{t-1}=[y_{t-1,1},\ldots,y_{t-1,\ell}]\) fixed,
apply the programmed preconditioner to \(Ay_{t-1,a}\)
for each column and subtract its output from \(y_{t-1,a}\),
\begin{equation}
    \label{eq:subspace_iteration_clean}
    \begin{aligned}
        u_{t,a}
        &=\left(I-B(W,\omega_{t,a})A\right)y_{t-1,a},
        \quad a=1,\ldots,\ell,\\
        U_t&=[u_{t,1},\ldots,u_{t,\ell}],  \\
        Y_t&=\operatorname{orth}(U_t),
    \end{aligned}
\end{equation}
where \(\omega_{t,a}\) denotes the fresh read collection
for this application, which averages \(S\) reads
from each of the same \(D\) programmed copies, and
each column uses a separate fresh read collection.
The products with \(A\), subtraction, and Euclidean
orthogonalization are digital.
Under the independent-read model,
\begin{equation}
    \label{eq:columnwise_probe_mean}
    \E_{\omega_{t,a}}\!\left[
        u_{t,a}\mid W,Y_{t-1}
    \right]
    =
    \overline E_Wy_{t-1,a}.
\end{equation}
Thus, the columns of \(U_t\) sample the mean-error
outputs, and \(Y_t\) provides an orthonormal basis
for their span. This basis supplies the probes for
the next sweep.
If fresh-read fluctuations are absent, then
\(U_t=\overline E_WY_{t-1}\), and the resulting candidate space
satisfies
$\range(Y_q)=\range(\overline E_W^{\,q}\Omega)$.
%
We next select correction directions within the
candidate space \(\range(Y_q)\).
Apply the same columnwise sampling
once more to \(Y_q\), using fresh reads to form
\(U_{q+1}\) and compute
\(H=Y_q^TU_{q+1}\in\R^{\ell\times\ell}\).
Thus, \(H\) estimates the mean-error action projected
onto the candidate space, i.e.,
\begin{equation}
        \E_\omega[H\mid W,Y_q]
        =Y_q^T\overline E_WY_q.
\end{equation}
Since \(H\) need not be symmetric, we use a real
Schur decomposition to obtain real orthonormal
selection directions:
$H=QTQ^T, Q^TQ=I$.
Here, \(T\) is real quasi-upper-triangular, with
\(1\times1\) blocks for real eigenvalues and
\(2\times2\) blocks for complex-conjugate pairs.
Reorder the Schur blocks by decreasing eigenvalue
modulus. The columns of \(Q\) are coordinates within
\(\range(Y_q)\), so the first \(k\) selected directions
give
$Z=Y_qQ(:,1{:}k)$.
If the cutoff after column \(k\) splits a
\(2\times2\) block, use \(Z=Y_qQ(:,1{:}k+1)\)
instead, retaining both associated columns and giving
an actual correction rank of \(k+1\).
The eigenvalue ordering targets modes with larger
moduli in the projected error operator; this selection
rule does not in general minimize
\eqref{eq:oracle_correction_space}.
The resulting basis remains fixed during the solve.
With no rank loss, the \(q\) candidate-space sweeps
and the additional operator-estimation sweep use
\(\ell(q+1)\) averaged analog vector applications.

Because the candidate space is built from \(U_t\), we examine how accurately
these columns represent
\(\overline E_WY_{t-1}\).
The following theorem gives conditional mean and
variance estimates for each implemented sampling
sweep. The programmed realization \(W\) is held
fixed throughout: the analysis concerns fresh read
noise in observing the mean-error action, not the
removal of its persistent discrepancies.
For \(X=[x_1,\ldots,x_\ell]\in\R^{N\times\ell}\),
write
$
    \norm{X}_{A,F}^2
    =
    \sum_{a=1}^{\ell}\norm{x_a}_A^2
    =
    \operatorname{tr}(X^TAX)$.
\begin{theorem}[Mean and variance of columnwise sampling]
\label{thm:programmed_sketch}
Assume the read model of
Section~\ref{sec:programmed_mean} for all columnwise
applications \(B(W,\omega_{t,a})\), each averaging \(S\) reads
from each of the same \(D\) programmed copies.
Conditional on \(W\), assume that fresh read
collections are independent across columns and
sweeps and independent of \(\Omega\).
Suppose that a finite constant \(\sigma_{\rm app}(W)\geq0\)
satisfies
\begin{equation}
    \label{eq:fluctuation_bound_for_sketch}
    \E_{\omega_d}\!\left[
        \norm{\Delta_d(\omega_d)Ax}_A^2
        \,\middle|\,W
    \right]
    \leq\sigma_{\rm app}^2(W)\norm{x}_A^2
\end{equation}
for every deterministic \(x\in\R^N\) and every
copy \(d\).
For each sampling sweep \(t=1,\ldots,q+1\),
including the final sweep used to form \(H\),
\begin{equation}
    \label{eq:implemented_probe_mean}
    \E_\omega[U_t\mid W,Y_{t-1}]
    =\overline E_WY_{t-1},
\end{equation}
and
\begin{equation}
    \label{eq:implemented_probe_variance}
    \begin{aligned}
        &\E_\omega\!\left[
            \norm{U_t-\overline E_WY_{t-1}}_{A,F}^2
            \,\middle|\,W,Y_{t-1}
        \right]\leq
        \frac{\sigma_{\rm app}^2(W)}{DS}
        \norm{Y_{t-1}}_{A,F}^2.
    \end{aligned}
\end{equation}
The expectations are over fresh setup reads with
\(W,Y_{t-1}\) fixed.
\end{theorem}

\begin{proof}
The mean identity \eqref{eq:implemented_probe_mean}
follows by collecting the columnwise identities
in \eqref{eq:columnwise_probe_mean}.
We now prove the variance bound
\eqref{eq:implemented_probe_variance}.
Write \(\Delta_{t,a}=B(W,\omega_{t,a})-\overline B_W\)
and set \(\nu=\sigma_{\rm app}^2(W)/(DS)\).
For a vector \(x\) determined before the current
application, conditioning on the preceding computation
fixes \(x\) and leaves the fresh read distribution
unchanged. Thus, \(\Delta_{t,a}Ax\) has zero
conditional mean.
Conditional independence eliminates cross terms
between the \(DS\) reads, and the \(S\) reads of
each copy have identical conditional distributions.
Consequently,
\begin{equation}
    \label{eq:averaged_fluctuation_bound}
    \begin{aligned}
        \E_\omega\!\left[
            \norm{\Delta_{t,a}Ax}_A^2
            \,\middle|\,W,x
        \right]&=
        \frac1{D^2S}\sum_{d=1}^D
        \E_{\omega_d}\!\left[
            \norm{\Delta_d(\omega_d)Ax}_A^2
            \,\middle|\,W,x
        \right]
        \leq\nu\norm{x}_A^2.
    \end{aligned}
\end{equation}
The same estimates hold conditional on the full
preceding computation. Applying them to each column
of \(Y_{t-1}\) proves
\eqref{eq:implemented_probe_variance}, since the
squared weighted Frobenius norm sums the squared
column norms.
\end{proof}

The single-sweep estimates explain the role of setup
read averaging. For fixed programmed copies, current
basis, and single-read distributions, increasing \(S\)
reduces the variance bound as \(1/S\) without changing
the mean-error action being sampled. Repeated reads
therefore help observe the programmed errors more
accurately.
%
To study how read noise accumulates across
repeated applications, we consider a separate
recurrence that omits orthonormalization.
This auxiliary
iteration is used only for analysis.
For \(a=1,\ldots,\ell\) and \(t=1,\ldots,q\),
\begin{equation}
    \label{eq:auxiliary_iteration}
    \begin{aligned}
        \widetilde y_{0,a}&=\Omega(:,a),\\
        \widetilde y_{t,a}
        &=(I-B(W,\omega_{t,a})A)\widetilde y_{t-1,a} \ .
    \end{aligned}
\end{equation}
Collect its columns in
\(\widetilde Y_t
=[\widetilde y_{t,1},\ldots,\widetilde y_{t,\ell}]\).
Without orthonormalization, we can identify the
conditional mean and bound the accumulated variance.

\begin{proposition}
\label{prop:unorthogonalized_sketch}
Under the assumptions of Theorem~\ref{thm:programmed_sketch},
suppose additionally that a finite constant \(\rho_W\geq0\)
satisfies
\begin{equation}
    \label{eq:mean_bound_for_sketch}
    \norm{\overline E_Wx}_A^2
    \leq\rho_W^2\norm{x}_A^2
\end{equation}
for every deterministic \(x\in\R^N\).
For the auxiliary iteration \eqref{eq:auxiliary_iteration}, we have
\begin{align}
    \label{eq:programmed_sketch_mean}
    \E_\omega[\widetilde Y_q\mid W,\Omega]
    &=\overline E_W^q\Omega, \\
    \label{eq:programmed_sketch_variance}
    \E_\omega\!\left[
            \norm{
                \widetilde Y_q-\overline E_W^q\Omega
            }_{A,F}^2
            \,\middle|\,W,\Omega
        \right] &\leq
        \left[
            \left(
                \rho_W^2+
                \frac{\sigma_{\rm app}^2(W)}{DS}
            \right)^q
            -\rho_W^{2q}
        \right]
        \norm{\Omega}_{A,F}^2.
\end{align}
The expectations are over fresh setup reads with
\(W,\Omega\) fixed.
\end{proposition}

\begin{proof}
Set \(\nu=\sigma_{\rm app}^2(W)/(DS)\).
The conditional centering and fluctuation bound
\eqref{eq:averaged_fluctuation_bound} apply to the
auxiliary columns as well. Conditional centering gives
\begin{equation}
    \E_\omega[
        \widetilde y_{t,a}\mid W,\Omega
    ]
    =
    \overline E_W
    \E_\omega[
        \widetilde y_{t-1,a}\mid W,\Omega
    ].
\end{equation}
Induction proves
\eqref{eq:programmed_sketch_mean} columnwise.
Define
\begin{equation}
    \epsilon_{t,a}
    =
    \E_\omega\!\left[
        \norm{
            \widetilde y_{t,a}
            -\overline E_W^t\Omega(:,a)
        }_A^2
        \,\middle|\,W,\Omega
    \right].
\end{equation}
After subtracting the column mean from the recurrence,
the cross term between the propagated previous error
and the current fluctuation vanishes conditional on
preceding reads. Using
\eqref{eq:mean_bound_for_sketch},
\eqref{eq:averaged_fluctuation_bound}, and the column
mean identity gives
\begin{equation}
    \begin{aligned}
        \epsilon_{t,a}
        &\leq
        \rho_W^2\epsilon_{t-1,a}
        +
        \nu\,\E_\omega\!\left[
            \norm{\widetilde y_{t-1,a}}_A^2
            \,\middle|\,W,\Omega
        \right]\\
        &=
        (\rho_W^2+\nu)\epsilon_{t-1,a}
        +
        \nu\norm{
            \overline E_W^{t-1}\Omega(:,a)
        }_A^2\\
        &\leq
        (\rho_W^2+\nu)\epsilon_{t-1,a}
        +
        \nu\rho_W^{2(t-1)}
        \norm{\Omega(:,a)}_A^2.
    \end{aligned}
\end{equation}
The equality decomposes the second moment into
variance and squared mean. Starting from
\(\epsilon_{0,a}=0\), induction yields
\begin{equation}
    \epsilon_{q,a}
    \leq
    \bigl[
        (\rho_W^2+\nu)^q-\rho_W^{2q}
    \bigr]\norm{\Omega(:,a)}_A^2.
\end{equation}
Summing over columns and substituting the definition
of \(\nu\) proves
\eqref{eq:programmed_sketch_variance}.
\end{proof}

Proposition~\ref{prop:unorthogonalized_sketch} bounds
read-noise accumulation over a chosen depth \(q\).
For fixed \(W,\Omega,q\) and single-read distributions,
the bound tends to zero as \(S\to\infty\).
These bounds apply to the auxiliary iterate
\(\widetilde Y_q\); they do not bound the error in
the final correction space obtained by orthonormalization
and Schur selection.

\subsection{Write-noise sensitivity and local regularization}
\label{sec:regularization}

Before programming the local inverses, we choose a
regularization shift \(\beta\). A larger shift can make the
programmed inverse less sensitive to write noise, but it also
moves the target away from the desired inverse \(A_j^{-1}\).
We therefore first identify shifts for which the write error
is small relative to the shifted target, and then examine the
approximation error introduced by the shift.

To obtain an explicit criterion, we specialize to centered
additive write noise and exclude fresh read noise, quantization,
and clipping. Fix an SPD block \(A_j\in\R^{n\times n}\),
where \(n=|\mathcal I_j|\), with orthonormal eigenvectors
satisfying
\begin{equation}
    \label{eq:local_block_eigenpairs}
    A_ju_l=\theta_lu_l,\qquad
    0<\theta_1\leq\cdots\leq\theta_n .
\end{equation}
For \(\beta\geq0\), define
\[
    G_\beta=(A_j+\beta I)^{-1},
    \qquad
    \mu=\frac{\norm{G_\beta}_{\max}}{w_{\rm r}},
\]
where \(\norm{X}_{\max}=\max_{i,l}|X_{il}|\) and
\(w_{\rm r}>0\) is the maximum intended weight magnitude.

Let \(\Phi_d\in\R^{n\times n}\) denote the additive
write error in copy \(d\). The programmed array stores
\(\mu^{-1}G_\beta+\Phi_d\), where the entries of
\(\Phi_d\) are independent with mean zero and variance
\(\sigma_{\rm w}^2\). The matrices are independent
across copies and remain fixed after programming.
Averaging the rescaled outputs of \(D\) copies gives
\[
    \widetilde G_\beta
    =
    G_\beta+\mu\overline\Phi_D,
    \qquad
    \overline\Phi_D
    =
    \frac1D\sum_{d=1}^D\Phi_d.
\]
Throughout this subsection, \(\E_\Phi\) denotes expectation
over programming realizations.

We require the root-mean-square write error to be smaller
than the intended output:
\begin{equation}
    \label{eq:write_noise_requirement}
    \sqrt{
        \E_\Phi
        \norm{(\widetilde G_\beta-G_\beta)r_j}_2^2
    }
    <
    \norm{G_\beta r_j}_2,
    \qquad r_j\neq0.
\end{equation}
A direction is called \emph{noise-dominated} when this
condition fails.

\begin{theorem}[Local write-noise criterion]
\label{thm:regularization}
Under the write-noise model above, define
\begin{equation}
    \label{eq:survival_fn}
    \tau_D
    =
    \frac{\sqrt D\,w_{\rm r}}
         {\sigma_{\rm w}\sqrt n},
    \qquad
    \varphi(\beta)
    =
    (\theta_n+\beta)\norm{G_\beta}_{\max}.
\end{equation}
For the unit eigendirection \(u_l\), the intended response
and root-mean-square write error are
\[
    \norm{G_\beta u_l}_2
    =
    \frac1{\theta_l+\beta},
    \qquad
    \sqrt{
        \E_\Phi
        \norm{(\widetilde G_\beta-G_\beta)u_l}_2^2
    }
    =
    \mu\sigma_{\rm w}\sqrt{\frac nD}.
\]
Thus, the write error is smaller than the intended response
in direction \(u_l\) if and only if
\begin{equation}
    \label{eq:survival_criterion}
    (\theta_l+\beta)\norm{G_\beta}_{\max}<\tau_D.
\end{equation}

Moreover, the worst relative root-mean-square write error
over deterministic nonzero inputs is
\begin{equation}
    \label{eq:sup_E}
    \sup_{r_j\neq0}
    \frac{
        \sqrt{
            \E_\Phi
            \norm{(\widetilde G_\beta-G_\beta)r_j}_2^2
        }
    }{
        \norm{G_\beta r_j}_2
    }
    =
    \frac{\varphi(\beta)}{\tau_D},
\end{equation}
with the maximum attained at \(r_j=u_n\).
Hence, \eqref{eq:write_noise_requirement} holds for every
deterministic nonzero input exactly when
\begin{equation}
    \label{eq:global_survival_criterion}
    \varphi(\beta)<\tau_D.
\end{equation}

The function \(\varphi\) is continuous and nonincreasing
on \([0,\infty)\), satisfies \(\varphi(\beta)\geq1\),
and converges to one as \(\beta\to\infty\).
A finite shift satisfying
\eqref{eq:global_survival_criterion} therefore exists
if and only if \(\tau_D>1\).
\end{theorem}

\begin{proof}
Independence and centering give, for every deterministic
\(r_j\),
\begin{equation}
    \label{eq:E_phi}
    \E_\Phi
    \norm{(\widetilde G_\beta-G_\beta)r_j}_2^2
    =
    \mu^2\sigma_{\rm w}^2\frac nD
    \norm{r_j}_2^2.
\end{equation}
For \(r_j=u_l\), this gives the stated root-mean-square
write error, while
\(G_\beta u_l=(\theta_l+\beta)^{-1}u_l\).
Combining the two expressions gives
\eqref{eq:survival_criterion}.

Since the smallest eigenvalue of \(G_\beta\) is
\(1/(\theta_n+\beta)\),
\[
    \norm{G_\beta r_j}_2
    \geq
    \frac{\norm{r_j}_2}{\theta_n+\beta},
\]
with equality for \(r_j=u_n\). Together with
\eqref{eq:E_phi}, this proves \eqref{eq:sup_E}.

Finally, because \(G_\beta\) is SPD,
its largest entry magnitude occurs on the diagonal, and
\[
    \varphi(\beta)
    =
    \max_i\sum_{l=1}^n
        (u_l)_i^2
        \frac{\theta_n+\beta}{\theta_l+\beta}.
\]
Each ratio is continuous, nonincreasing, at least one,
and converges to one. Since the weights
\((u_l)_i^2\) are nonnegative and sum to one,
the same properties hold for \(\varphi\).
The feasibility statement follows from
\(\varphi(\beta)\to1\).
\end{proof}

Theorem~\ref{thm:regularization} compares two quantities:
the intended response \(1/(\theta_l+\beta)\) and the
write-error level
\(\mu\sigma_{\rm w}\sqrt{n/D}\).
The criterion is satisfied in direction \(u_l\) when the
former lies above the latter. Increasing \(D\) lowers the
write-error level as \(D^{-1/2}\). Changing \(\beta\)
changes both quantities through \(G_\beta\) and the
programming scale \(\mu\).

This criterion measures how accurately the hardware
implements \(G_\beta\). A separate issue is how accurately
\(G_\beta\) approximates \(A_j^{-1}\).
To examine this, consider instead the input
\(r_j=A_ju_l=\theta_lu_l\), for which the exact local
correction is \(u_l\). The shifted inverse leaves the
deterministic error
\begin{equation}
    \label{eq:regularization_error}
    u_l-G_\beta r_j
    =
    \frac{\beta}{\theta_l+\beta}u_l.
\end{equation}
This error increases with \(\beta\). Theorem~\ref{thm:regularization} shows that increasing
\(\beta\) does not worsen the relative write-noise
criterion. This does not, however, imply that a larger
shift gives a better preconditioner.

We therefore use \eqref{eq:global_survival_criterion}
as a feasibility condition and choose a feasible value
near its lower boundary. If \(\tau_D\leq1\), no finite
shift satisfies the criterion. For \(\tau_D>1\), define
\begin{equation}
    \label{eq:beta_min}
    \beta_{\min}
    :=
    \inf\{
        \beta\geq0:
        \varphi(\beta)<\tau_D
    \}.
\end{equation}
If \(\varphi(0)<\tau_D\), we set
\(\beta_{\text{safe},j}=0\).
Otherwise, continuity gives
\(\varphi(\beta_{\min})=\tau_D\), while monotonicity
and the definition of the infimum imply that every
\(\beta>\beta_{\min}\) satisfies the strict criterion.

Numerically, we bracket \(\beta_{\min}\) by an infeasible
lower endpoint and a feasible upper endpoint and apply
bisection. For a prescribed tolerance
\(\varepsilon_\beta\), the final feasible upper endpoint
is used as \(\beta_{\text{safe},j}\).
If a common shift is required for all blocks, we take
\[
    \beta=\max_j\beta_{\text{safe},j}.
\]
This gives a feasible reference shift under the stated
write-noise model; it is not an iteration-optimal choice
for FGMRES.

\begin{figure}[t]
    \centering
    \includegraphics[width=0.95\linewidth]{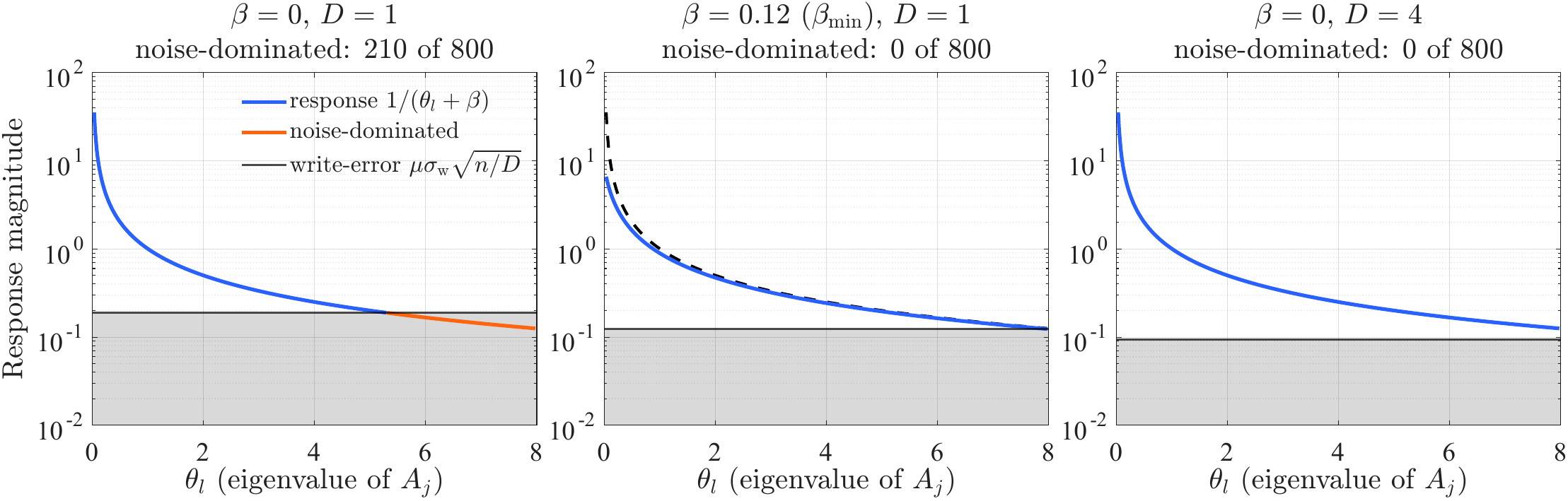}
    \caption{Local write-noise comparison for the first
    block of a 2D Poisson problem on a
    \(40\times40\) grid, with \(p=2\), \(n=800\), and
    noise factor \(f_n=1.2\).
    The curves show the intended responses
    \(1/(\theta_l+\beta)\), and the black horizontal line
    shows the root-mean-square write-error level
    \(\mu\sigma_{\rm w}\sqrt{n/D}\).
    Red segments indicate noise-dominated directions.
    Left: \(\beta=0\), \(D=1\), with \(210\) such
    directions.
    Middle: \(\beta\approx0.12\), \(D=1\), near the
    predicted threshold; the dashed curve shows the
    unregularized response.
    Right: \(\beta=0\), \(D=4\); copy averaging halves
    the write-error level without changing the target
    responses.}
    \label{fig:dynrange}
\end{figure}

Figure~\ref{fig:dynrange} visualizes
\eqref{eq:survival_criterion} directly.
The response curve shows
\(1/(\theta_l+\beta)\), while the black horizontal line
shows the root-mean-square write-error level
\(\mu\sigma_{\rm w}\sqrt{n/D}\).
A direction is noise-dominated when its response lies
at or below the black line.

\section{Numerical Experiments}
\label{sec:experiments}

The experiments examine whether the design choices motivated
by Section~\ref{sec:theory} improve the simulated preconditioner.
We combine error diagnostics, which examine quantities appearing
in the analysis, with solver comparisons, which measure their
practical consequences for FGMRES.

On a \(40\times40\) Poisson problem (\(N=1600\), discretized using the standard 5-point finite difference stencil),
Section~\ref{sec:exp-lowrank} tests whether programming changes
the correction directions and whether learning from the
programmed copies improves on construction from digital
block Jacobi. Section~\ref{sec:exp-biasvar} examines the
programmed-mean error, then relates a combined mean-plus-variance score to iteration counts.
Section~\ref{sec:exp-large} tests correction and averaging
on larger Poisson and interface problems.
Finally, Section~\ref{sec:exp-fem} compares the local
write-noise threshold with the regularization shifts
that improve solver performance.

All experiments use \textsc{Matlab} with 64-bit arithmetic
on all eight cores of an AMD Ryzen AI 7 PRO 350 CPU.
Analog operations are simulated using a \textsc{Matlab}
implementation of the simulator described in \cite{analogsim}.
The default parameters in Table~\ref{tab:params} follow
previous analog preconditioning studies
\cite{kalantzis2021solving,kalantzis2023fgmres,shah2026hybrid}.
The noise factor \(f_n\) multiplies the standard deviations
of write, input, and output noise, and of the device weight range variation.

Unless stated otherwise, we use \(M_Z(B)\) from
\eqref{eq:low_rank_corrected_action} with
\(\beta=0\) and correction rank \(k=32\).
The basis is constructed using \(q=10\) columnwise
sampling sweeps of \eqref{eq:subspace_iteration_clean}
with probe width \(\ell=3k\).
Setup and solve use the same \(D\) programmed copies
and the same averaging: each vector application uses
\(S\) reads per copy.

Unrestarted FGMRES uses a zero initial guess, a relative residual tolerance
of \(10^{-8}\), and a limit of \(1000\) iterations.
Unless stated otherwise, results report means over
\(10\) independent programming realizations, each with
one random unit-norm right-hand side. Shaded bands
extend one standard deviation above and below the mean
across these trials.

\begin{table}[t]
    \centering
    \caption{Default simulator parameters. For each local
    target \(G_\beta=(A_j+\beta I)^{-1}\), scaling places
    the entries of \(\mu^{-1}G_\beta\) in
    \([-w_{\rm r},w_{\rm r}]\).}
    \label{tab:params}
    \begin{tabular}{l c l}
        \toprule
        Parameter & Symbol & Value \\
        \midrule
        Noise factor
        & \(f_n\) & \(1.0\) (default) \\
        Write-noise standard deviation
        & \(\sigma_{\rm w}\) & \((5\times10^{-3})f_n\) \\
        Input-noise standard deviation
        & \(\sigma_I\) & \((1\times10^{-2})f_n\) \\
        Output-noise standard deviation
        & \(\sigma_O\) & \((1\times10^{-2})f_n\) \\
        Device weight range variation
        & \(\sigma_{\rm r}\) & \((2\times10^{-1})f_n\) \\
        \midrule
        Maximum weight magnitude
        & \(w_{\rm r}\) & \(0.6\) \\
        Matrix scaling factor
        & \(\mu\) & \(\norm{G_\beta}_{\max}/w_{\rm r}\) \\
        Number of DAC bits
        & \(\ell_{\rm DAC}\) & \(7\) \\
        Number of ADC bits
        & \(\ell_{\rm ADC}\) & \(9\) \\
        \bottomrule
    \end{tabular}
\end{table}

\subsection{Evaluation of the low-rank correction preconditioner}
\label{sec:exp-lowrank}

We first test the central construction choice: whether
the correction space should be learned from the programmed
analog action rather than digital block Jacobi.
On the \(40\times40\) Poisson problem, we examine how
programming changes the errors left by the base preconditioner,
whether the learned space responds to those changes,
and whether this response improves FGMRES performance.

By \eqref{eq:corrected_error_identity},
the error left by the analog base acts as the target of the
digital correction. We first examine this error before
adding the correction, excluding fresh input/output noise
to focus on the programmed action.
Let \(Av_i=\lambda_i v_i\), with Euclidean-orthonormal
eigenvectors ordered by increasing eigenvalue.
For either the digital or programmed analog operator \(B\),
define
\begin{equation}
    \rho_i
    =
    \frac{\norm{v_i-BAv_i}_A}{\norm{v_i}_A}.
\end{equation}
For input \(Av_i\), the exact correction is \(v_i\),
so \(\rho_i\) measures the relative approximation error.
A value \(\rho_i\geq1\) means that the error is at least
as large as the exact correction.

\begin{figure}[tbhp]
    \centering
    \includegraphics[width=0.95\textwidth]{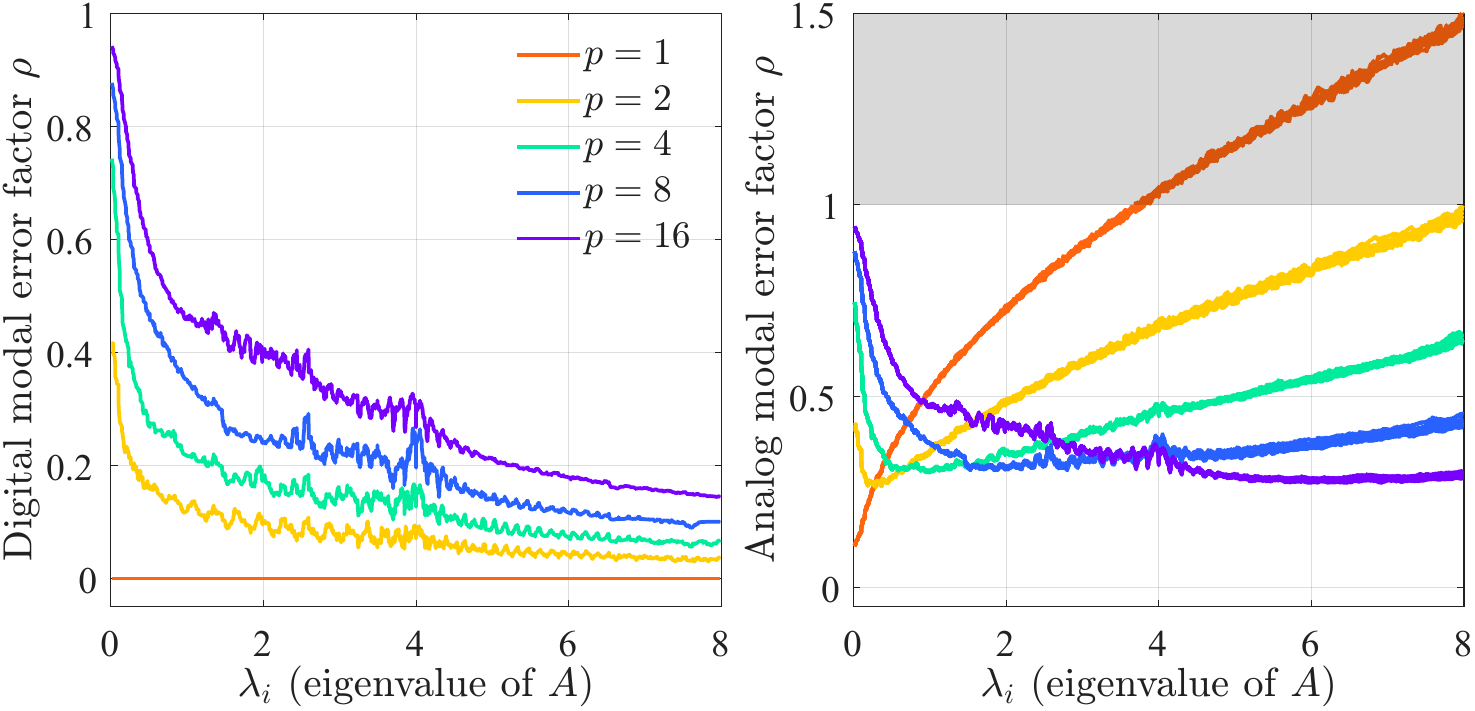}
    \caption{Modal error factors \(\rho_i\) for digital
    (left) and programmed analog (right) block Jacobi
    on the \(40\times40\) Poisson grid.
    The analog diagnostic includes write noise and
    device-dependent weight ranges at the default
    noise level, but excludes fresh input/output noise and quantization.
    Ten programming realizations are overlaid, with
    curves smoothed over \(16\) consecutive modes.
    The region \(\rho_i\geq1\) is shaded.
    Digital block Jacobi is exact for \(p=1\).}
    \label{fig:modal}
\end{figure}

Figure~\ref{fig:modal} shows that digital block Jacobi
mainly leaves low-eigenvalue directions poorly corrected:
at the smallest eigenvalue, \(\rho_i\) reaches \(0.99\)
for \(p=16\) (unsmoothed). Programming also introduces large errors
in high-eigenvalue directions. For \(p=1\), the analog
curves cross \(\rho_i=1\) near \(\lambda_i\approx3.8\),
with \(\rho_i\geq1\) for \(56\%\) of the modes.
Thus, a space designed only for the digital approximation
may miss errors introduced by programming.

\begin{figure}[tbhp]
    \centering
    \includegraphics[width=0.95\textwidth]{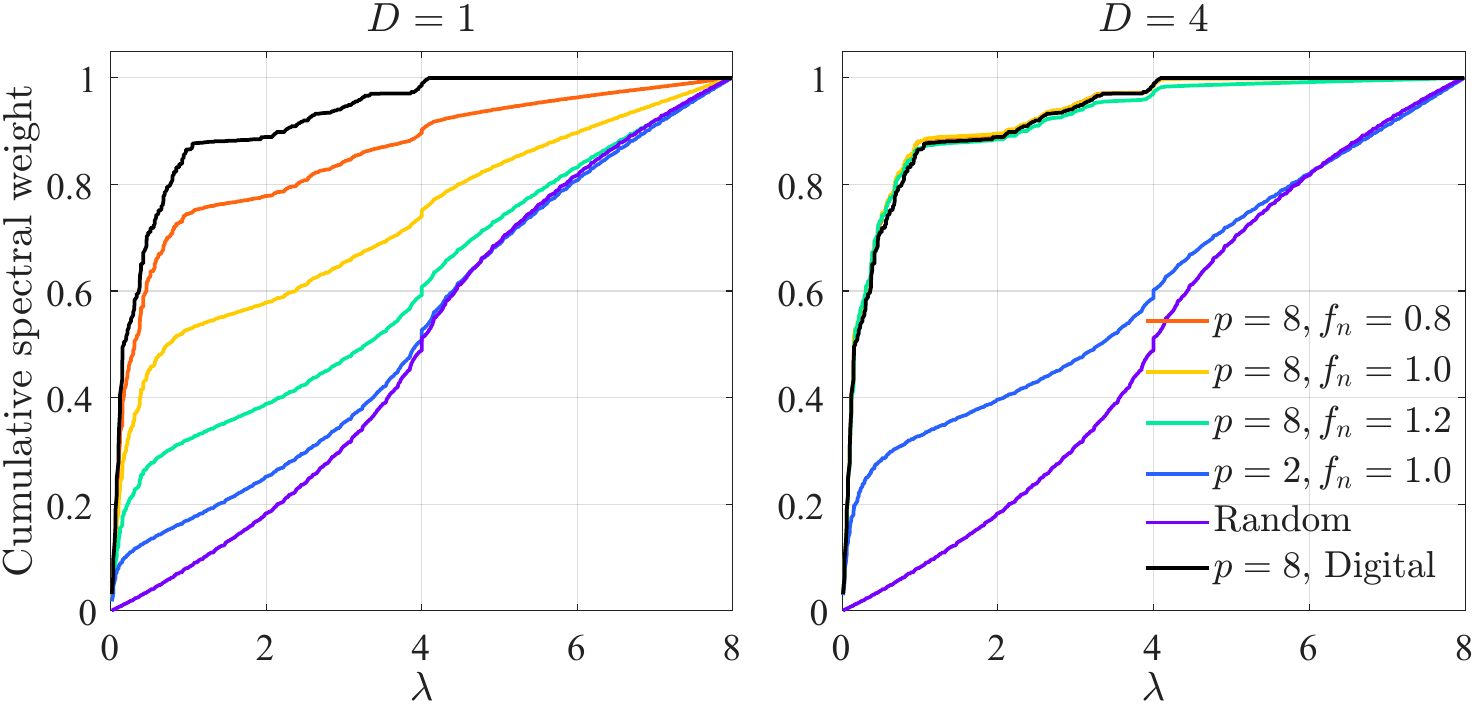}
    \caption{Cumulative spectral weight of rank-\(32\)
    correction spaces on the \(40\times40\) Poisson
    grid, using \(q=20\) sampling sweeps and probe width \(\ell=96\).
    Black: digital block Jacobi with \(p=8\).
    Colored curves: analog-learned spaces at the
    indicated block counts and noise levels, averaged
    over \(10\) programming realizations.
    Purple: random-subspace reference.
    Left: \(D=S=1\). Right: \(D=4\), \(S=1\).}
    \label{fig:energy}
\end{figure}

We next ask whether the space produced by probing responds
to these changes. Under the analysis model,
Theorem~\ref{thm:programmed_sketch} identifies the mean-error action sampled by each sweep.
Figure~\ref{fig:energy} examines the space obtained
after orthonormalization and Schur selection.
For a Euclidean-orthonormal basis \(Z\) with \(k\)
columns, its cumulative spectral weight is
\begin{equation}
    \frac1k\sum_{\lambda_i\leq\lambda}
        \norm{Z^Tv_i}_2^2.
\end{equation}
This is the fraction of the basis's total squared
Euclidean norm in eigenmodes below the cutoff \(\lambda\).
A large value at a small cutoff indicates concentration
in low-eigenvalue directions.

In Figure~\ref{fig:energy}, the digital space places
\(87\%\) of its spectral weight in \(\lambda\leq1\),
although these modes contain only \(8\%\) of the
eigenvectors and receive \(8\%\) expected weight from
an isotropically random subspace.
With \(p=8,D=1\), the analog-learned weights are
\(74\%\), \(53\%\), and \(32\%\) at
\(f_n=0.8,1.0,1.2\), respectively.
The learned space therefore becomes less concentrated
in low-eigenvalue directions as noise increases.
With \(D=4\), these curves approach the digital
reference. In the displayed \(p=2\) configuration,
increasing \(D\) from one to four raises the
low-eigenvalue weight from \(17\%\) to \(33\%\).
These measurements show that the final space changes
with programming noise and copy averaging.

\begin{figure}[tbhp]
    \centering
    \includegraphics[width=0.95\textwidth]{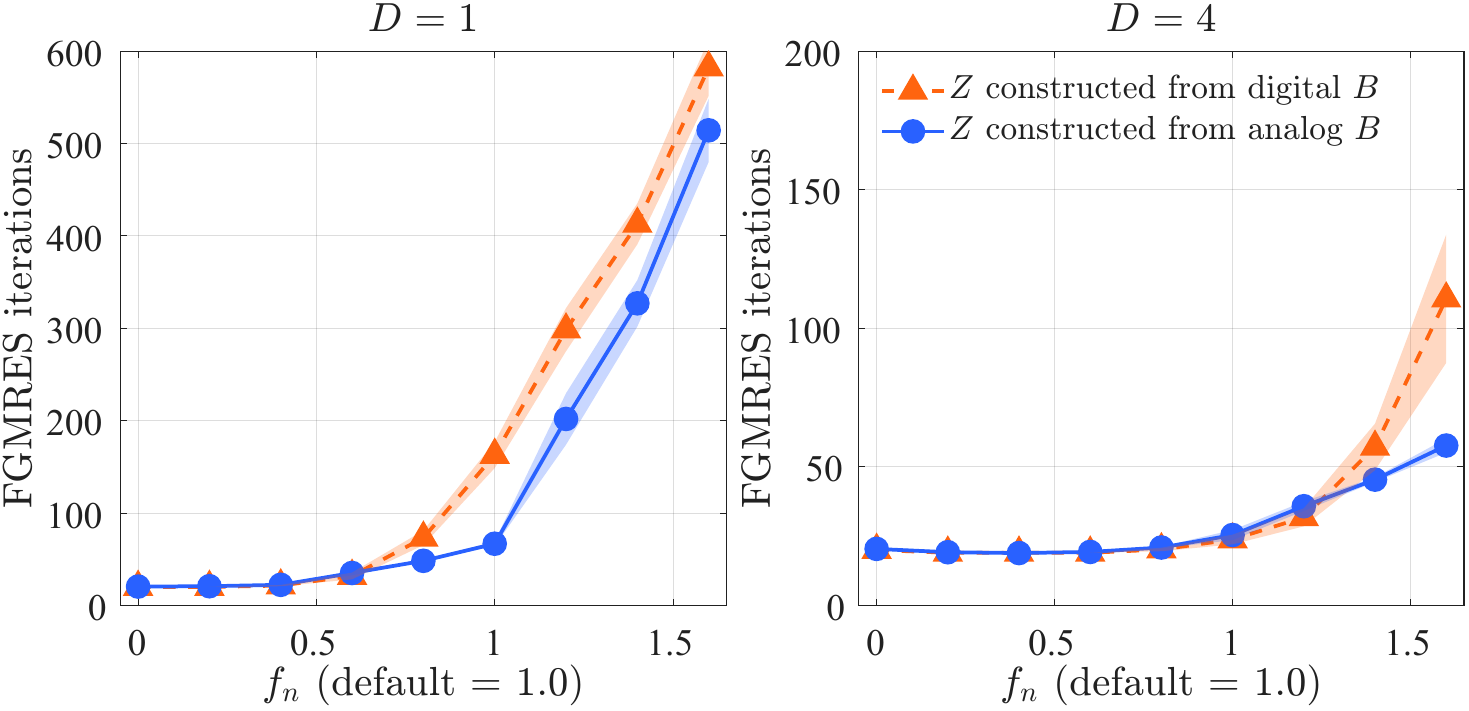}
    \caption{FGMRES iterations for correction spaces
    constructed from digital block Jacobi (red) or
    the programmed analog copies (blue), on the
    \(40\times40\) Poisson grid with \(p=4\),
    \(k=32\), and \(q=10\).
    Both methods use the same programmed analog
    base during the solve; analog setup uses the
    same copies and averaging as the solve.
    Left: \(D=S=1\). Right: \(D=4\), \(S=1\).}
    \label{fig:zsource}
\end{figure}

A change in the learned space does not by itself
establish a solver benefit. Figure~\ref{fig:zsource}
therefore compares spaces constructed from digital
block Jacobi and the programmed analog action while
holding the correction rank and programmed analog
base fixed during the solve.
Only the operator used to construct \(Z\) changes.
We use \(p=4\) and vary \(f_n\) from zero to \(1.6\).

For \(D=1\), Figure~\ref{fig:zsource} shows similar
iteration counts up to \(f_n=0.6\).
At the default noise level, however, the analog-learned
space reduces the mean count from \(162.9\pm14.2\)
to \(66.8\pm2.5\).
For \(D=4\), the two constructions perform similarly
up to \(f_n=1.2\); at \(f_n=1.6\), the counts are
\(110.7\pm23.2\) for the digital space and
\(57.7\pm2.5\) for the analog-learned space.
In these tests, learning from the programmed action
is most beneficial when programming perturbations
are significant and fewer copies are averaged.

\subsection{Mean--variance diagnostics and solver performance}
\label{sec:exp-biasvar}

\begin{figure}[tbhp]
    \centering
    \includegraphics[width=0.95\textwidth]{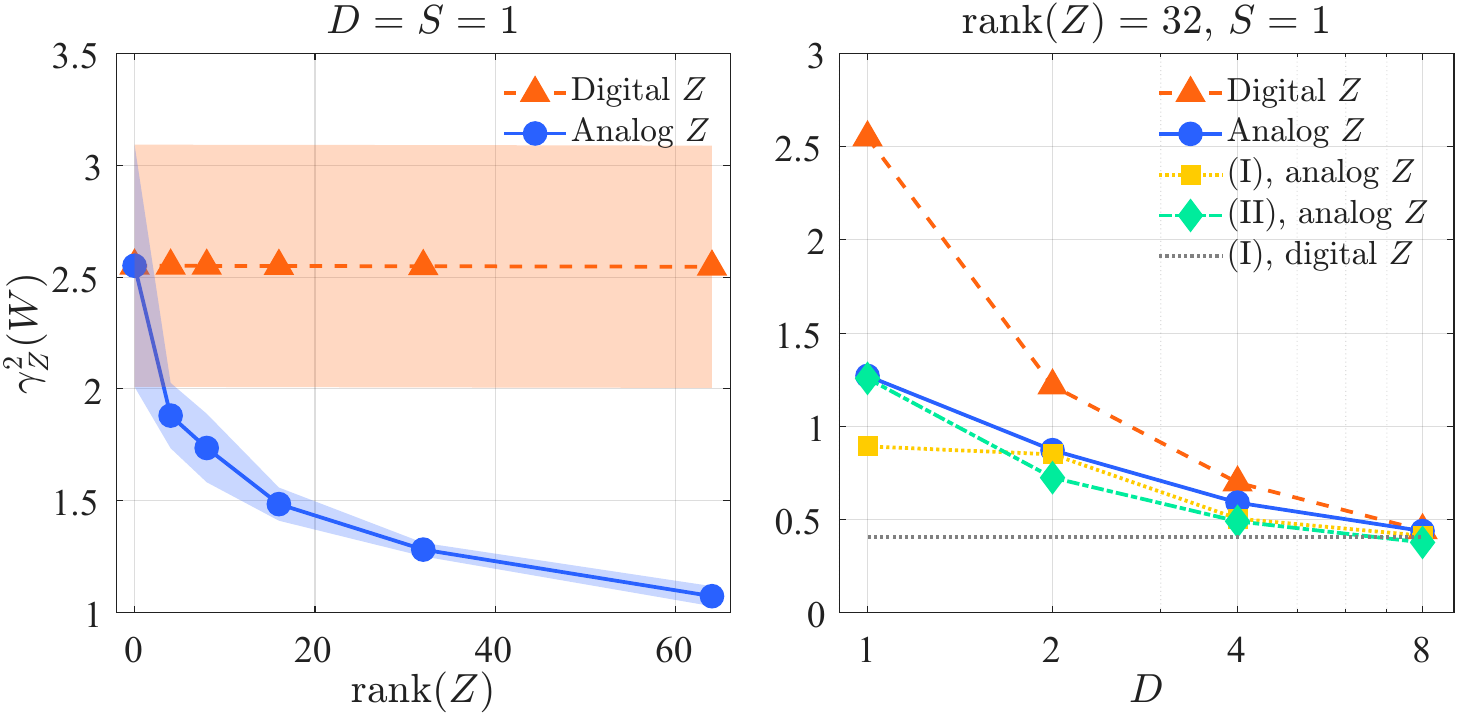}
    \caption{Offline programmed-mean error diagnostics on the \(40\times40\) Poisson
    problem with \(p=4\), \(\beta=0\), and quantization
    disabled.
    Left: \(\gamma_Z^2(W)\) versus correction rank
    for \(D=S=1\), using digital (red) or
    analog-learned (blue) spaces.
    Right: \(\gamma_Z^2(W)\) versus \(D\) for both
    correction spaces at rank \(32\) and \(S=1\),
    together with component (I) for both spaces and
    component (II) for the analog-learned space.
    Results are means over \(10\) programming realizations;
    bands (left only) show \(\pm\) one standard deviation.}
    \label{fig:biasvar}
\end{figure}

\begin{figure}[tbhp]
    \centering
    \includegraphics[width=0.50\textwidth]{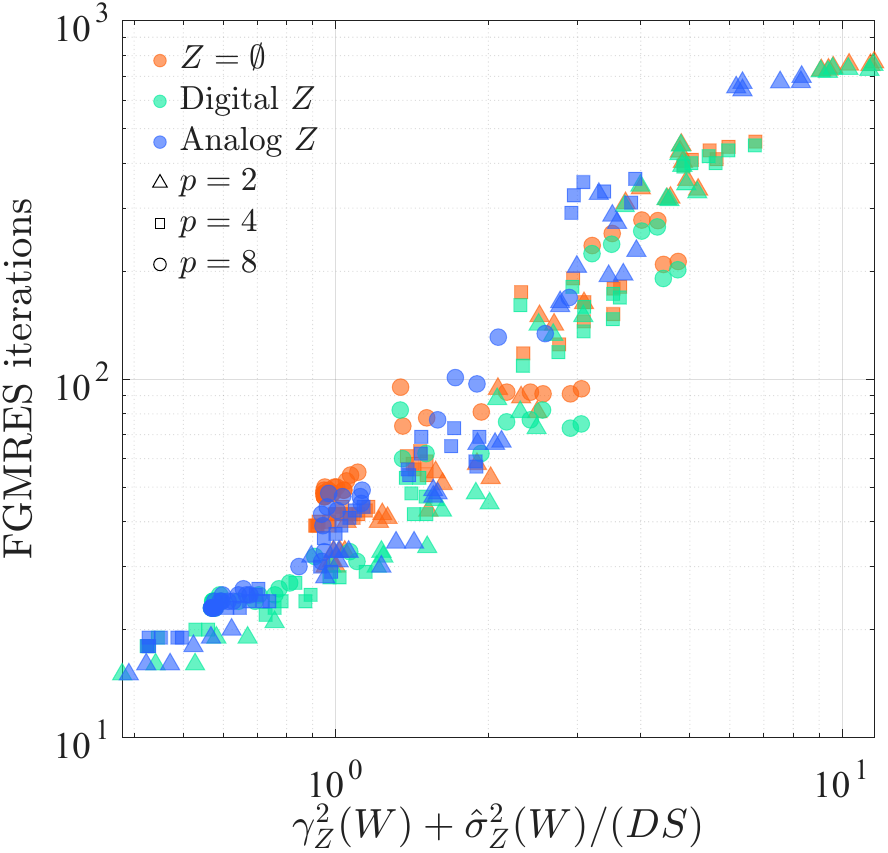}
    \caption{FGMRES iterations versus the offline
    diagnostic score
    \(\gamma_Z^2(W)+\widehat\sigma_Z^2(W)/(DS)\)
    on the \(40\times40\) Poisson problem.
    Each of the \(324\) markers represents one
    configuration and programming realization,
    with the median iteration count over three
    right-hand sides.
    Shape identifies \(p\), and color identifies
    the correction choice.
    Corrected cases use rank \(32\);
    quantization is disabled.}
    \label{fig:spearman}
\end{figure}

We now examine which error contributions help explain
the solver improvements observed in the preceding subsection.
Figure~\ref{fig:biasvar} measures the programmed-mean
error.
Figure~\ref{fig:spearman} then tests whether the 
combined diagnostic score is associated with FGMRES
iteration counts.

These are offline diagnostics on the \(40\times40\)
Poisson problem (\(N=1600\)). The correction-space
construction in Section~\ref{sec:learn_programmed_space}
requires only vector applications of the programmed
preconditioner; it does not require the entries of
\(\overline B_W\). For this small simulator test,
however, the effective programmed matrices can be
formed explicitly to evaluate the completed spaces.
This additional matrix access is used for diagnosis,
not as an input to the probing construction.

DAC/ADC quantization is disabled, and the remaining
input/output perturbations are additive and zero mean.
For the diagnostic calculations, we take
\(\overline B_W\) to be the average of the effective
programmed block Jacobi matrices over the \(D\)
copies within a fixed realization \(W\).
This identification requires the remaining read
perturbations to have zero conditional mean around
those matrix actions. We then explicitly form
\(\overline E_W=I-\overline B_WA\).

For Figure~\ref{fig:biasvar}, we fix \(p=4\) and
\(\beta=0\). Given a completed basis \(Z\), define
\begin{equation}
    \gamma_Z^2(W)
    :=
    \norm{(I-P_Z)\overline E_W}_A^2
    =
    \norm{
        A^{1/2}(I-P_Z)\overline E_WA^{-1/2}
    }_2^2.
\end{equation}
This measures the largest squared error of the corrected
mean action over \(\norm{v}_A=1\).
It is the square of the objective in
\eqref{eq:oracle_correction_space} and the maximum
of the first term in
\eqref{eq:one_step_bias_variance_clean}.
We compute it as the squared largest singular value
of the explicitly formed transformed operator,
taking \(P_Z=0\) when no correction is used.
This matrix calculation evaluates the space after
construction.

Figure~\ref{fig:biasvar} (left) shows that the
analog-learned space reduces this mean-error objective
more effectively. The digital space leaves
\(\gamma_Z^2(W)\) near \(2.55\) as the rank increases
from zero to \(64\), whereas the analog-learned
space reduces it to \(1.28\pm0.03\) at rank \(32\)
and \(1.07\pm0.04\) at rank \(64\).

The curves labeled ``Digital Z'' and ``Analog Z''
in Figure~\ref{fig:biasvar} (right) show the
combined programmed-mean error \(\gamma_Z^2(W)\)
for the two correction spaces.
The additional curves show selected component
norms from \eqref{eq:mean_error_split}.
Using \(\overline\Xi_D=\overline B_W-M_{\beta,p}\),
these quantities are
\begin{equation}
    \begin{aligned}
        \text{(I)}\quad&
        \norm{(I-P_Z)(I-M_{\beta,p}A)}_A^2,
        &&\text{digital approximation},\\
        \text{(II)}\quad&
        \norm{(I-P_Z)\overline\Xi_DA}_A^2,
        &&\text{programming discrepancy}.
    \end{aligned}
\end{equation}
Component (I) measures the approximation error
of ideal block Jacobi remaining after correction;
component (II) measures the remaining persistent
programming discrepancy. The figure shows
component (I) for both correction spaces and
component (II) only for the analog-learned space.

At \(D=1\), component (I) is \(0.41\) for the
digital space and \(0.89\) for the analog-learned
space. Thus, the digital space more effectively
corrects the ideal block Jacobi error in this
comparison. This does not determine which space
has the smaller combined programmed-mean error:
\(\gamma_Z^2(W)\) includes both error sources,
and the component norms are not additive because
their operators are subtracted before taking the
squared norm. At \(D=8\), the two combined-error
curves approach \(\gamma_Z^2(W)\approx0.44\).

Theorem~\ref{thm:one_step_bias_variance_clean}
gives an identity for the error of one corrected
application, but does not relate that error to
FGMRES iteration counts. We therefore examine
this relationship empirically.
Motivated by the two contributions in the identity,
we define the offline diagnostic score
\begin{equation}
    \gamma_Z^2(W)+\frac{\widehat\sigma_Z^2(W)}{DS}
\end{equation}
and test whether smaller scores are associated
with fewer iterations across the configurations
in Figure~\ref{fig:spearman}.
The first term is computed from the explicit
programmed-mean matrix as above.
To obtain \(\widehat\sigma_Z^2(W)\), we estimate
the projected single-read variance for each test
direction, average it over the \(D\) copies,
and take the largest of these directionwise estimates.
The test set contains \(31\) random \(A\)-unit vectors and three directions
associated with the largest corrected programmed-mean errors, one for each
correction choice; the same \(34\) directions are used for all three choices.
These directions are used only for the offline
variance diagnostic, not to construct \(Z\).

The parameter study uses \((D,S)\in\{(1,1),(1,4),(4,1),(2,2)\}\), \(p\in\{2,4,8\}\), and \(f_n\in\{0.6,1.0,1.4\}\).
We compare no correction, a digitally constructed
space, and an analog-learned space; corrected cases
use rank \(32\).
Three programming realizations per configuration
give \(324\) cases, each reporting the median
iteration count over three right-hand sides.
Quantization remains disabled.

Across the configurations in Figure~\ref{fig:spearman},
smaller scores are generally associated with fewer
iterations. The Spearman rank correlation is \(0.96\),
with iteration counts ranging from \(15\) to \(766\).
The median plotted counts are \(39\), \(48\), and
\(58.5\) for analog-learned, digital, and uncorrected
cases, respectively.
The score therefore helps interpret solver performance
in these tests.

\subsection{Larger 2D Poisson and domain decomposition}
\label{sec:exp-large}
\begin{figure}[tbhp]
    \centering
    \includegraphics[width=0.5\textwidth]{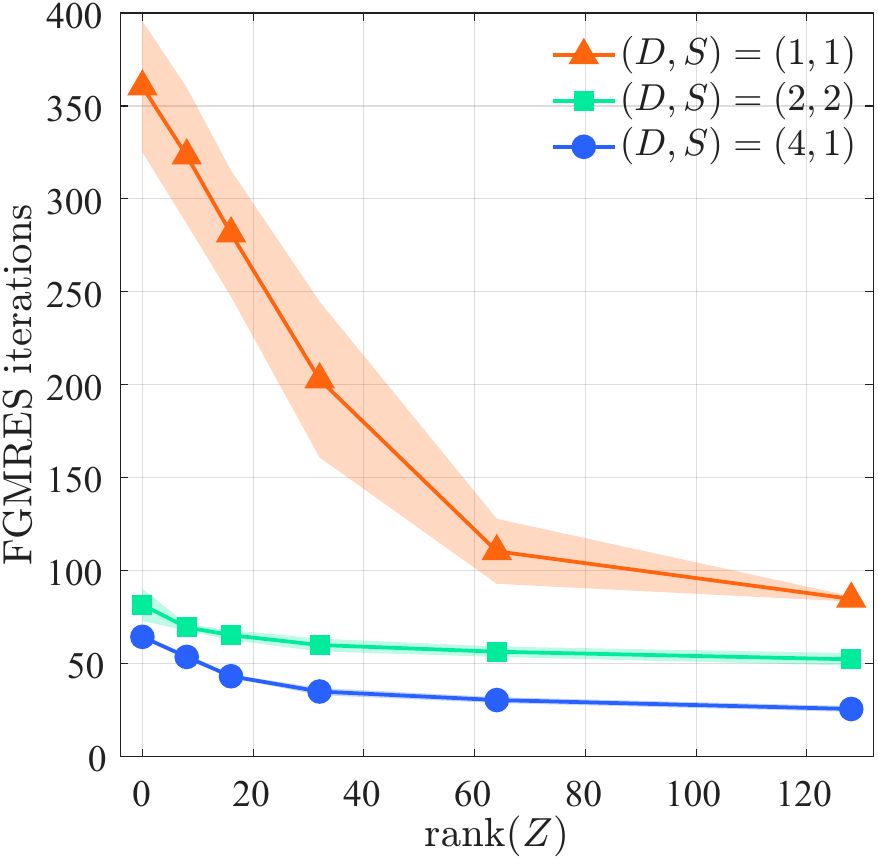}
    \caption{FGMRES iterations versus correction rank
    on the \(64\times64\) Poisson grid with \(p=8\),
    the default noise profile, and quantization enabled.
    Setup and solve use the same programmed copies
    and averaging.
    Markers show means over \(10\) trials;
    shaded bands indicate one standard deviation.}
    \label{fig:poisson64}
\end{figure}

We next examine whether correction and copy averaging
remain useful beyond the small diagnostic problem.
The first test uses a \(64\times64\) Poisson grid
(\(N=4096\)), with \(p=8\), the default noise profile,
and DAC/ADC quantization enabled.
Figure~\ref{fig:poisson64} varies the correction rank
from zero to \(128\) for
\((D,S)=(1,1),(2,2),(4,1)\).
The latter two configurations have the same number
of array reads per application, allowing us to
compare more independently programmed copies with
more repeated reads.
Within each trial, the configurations use subsets
of the same collection of programmed arrays.

For \(D=S=1\), Figure~\ref{fig:poisson64} shows that
the correction reduces the iteration count from
\(360.4\pm35.3\) without correction to
\(110.5\pm17.5\) at rank \(64\) and
\(84.9\pm1.4\) at rank \(128\).
At rank \(32\), four copies with one read each
require \(35.1\pm1.9\) iterations, compared with
\(60.1\pm3.4\) for two copies with two reads each.
This advantage is consistent with the distinction
in \eqref{eq:mean_error_copy_averaging}: copy
averaging can reduce programming variability,
whereas repeated reads leave the programmed
mean unchanged.
The comparison holds array-read counts fixed,
not hardware cost, since the configurations
use different numbers of arrays and sequential
read rounds.

The second test examines whether probing can construct
a useful correction when the analog base is built
from a matrix different from that of the system being solved.
We use the interface Schur complement of a
\(256\times256\) Poisson problem partitioned into
\(8\times8\) subdomains:
\begin{equation}
    \mathcal S
    =
    C-A_{\Gamma I}A_{II}^{-1}A_{I\Gamma},
    \qquad
    C=A_{\Gamma\Gamma}.
\end{equation}
Here, \(I\) denotes the interior unknowns,
\(\Gamma\) the interface, and \(\mathcal S\)
has dimension \(3535\times 3535\).

The analog base applies block Jacobi constructed
from \(C\), with \(p=7\) and \(D=S=1\).
Each probe nevertheless uses \(v-B\mathcal Sv\):
for input \(\mathcal Sv\), the exact correction
is \(v\). The sampled error therefore corresponds
to the Schur system, not to \(C\) alone.
The digital reference applies \(C^{-1}\) exactly,
without a low-rank correction.

\begin{table}[t]
    \centering
    \caption{FGMRES iterations for the interface Schur
    complement \(\mathcal S\).
    The analog base uses block Jacobi on \(C\), with
    \(p=7\) and \(D=S=1\); the correction space is
    probed using \(\mathcal S\).
    Means and standard deviations are computed over
    \(10\) trials, each with one right-hand side and,
    for the analog cases, an independent programming
    realization.}
    \label{tab:schur}
    \begin{tabular}{l c r r}
        \toprule
        & & \multicolumn{2}{c}{FGMRES iterations} \\
        \cmidrule(lr){3-4}
        Preconditioner & \(\rank(Z)\)
        & Mean & Std.\ dev. \\
        \midrule
        Digital \(C^{-1}\)
        & --- & 94.0 & 0.47 \\
        \midrule
        Analog block Jacobi
        & 0   & 142.1 & 2.28 \\
        With learned correction
        & 8   & 129.3 & 3.95 \\
        & 16  & 112.7 & 4.06 \\
        & 32  & 82.1  & 1.85 \\
        & 64  & 47.5  & 0.85 \\
        & 128 & 34.6  & 0.97 \\
        \bottomrule
    \end{tabular}
\end{table}

Table~\ref{tab:schur} shows that a rank-\(32\)
correction lowers the analog iteration count from
\(142.1\pm2.3\) to \(82.1\pm1.9\), below the
uncorrected digital reference of \(94.0\pm0.5\).
Increasing the rank to \(128\) reduces the count
to \(34.6\pm1.0\).
These results demonstrate the benefit of learning
a correction for the Schur-system error left
by the \(C\)-based analog approximation.

\subsection{SuiteSparse finite element stiffness matrix}
\label{sec:exp-fem}

Theorem~\ref{thm:regularization} identifies a threshold
for relative write error under its local noise model.
We now examine how this reference threshold relates
to the shifts that improve the complete simulated
preconditioner.
In particular, we compare the onset of successful
solves and the best tested shifts with
\(\beta_{\min}(D)\), and examine whether increasing
\(\beta\) further remains beneficial.

We use the SuiteSparse \cite{suitesparse} stiffness matrix
\texttt{bcsstk21} (\(N=3600\)), with symmetric Jacobi
scaling \(A\leftarrow\Lambda^{-1/2}A\Lambda^{-1/2}\),
where \(\Lambda=\operatorname{diag}(A)\).
We set \(p=8\), the correction rank to \(72\), and \(S=1\),
and vary \(\beta\) from \(0.01\) to \(2\) for
\(D=1,2\).
The correction space is learned from the programmed
copies used in the solve.
The reference is unregularized digital block Jacobi
with a digitally constructed space of the same rank.

\begin{figure}[tbhp]
    \centering
    \includegraphics[width=0.5\textwidth]{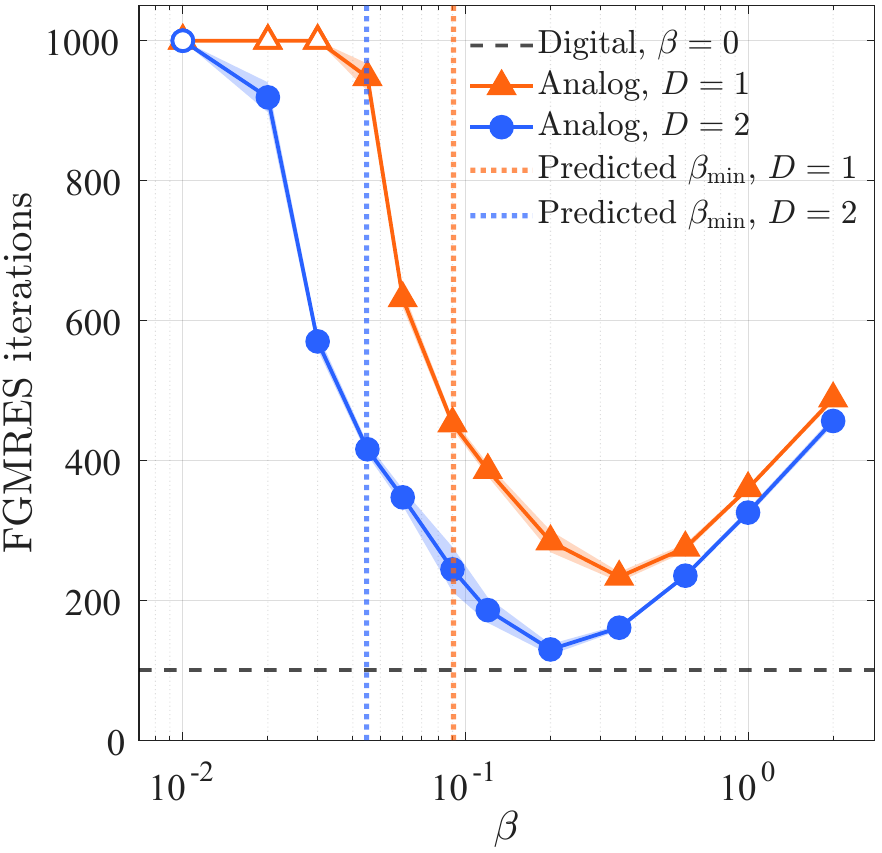}
    \caption{FGMRES iterations versus regularization
    shift for Jacobi-scaled \texttt{bcsstk21}, with
    \(p=8\), correction rank \(72\), \(S=1\), and
    \(D=1,2\).
    The dashed line is the unregularized digital
    reference; dotted lines mark the predicted local
    write-noise thresholds \(\beta_{\min}(D)\).
    Hollow markers indicate that at least one of
    the \(10\) runs did not converge within
    \(1000\) iterations.}
    \label{fig:bcsstk21}
\end{figure}

For each \(D\), let \(\beta_{\min}(D)\) denote
the largest blockwise threshold defined in
\eqref{eq:beta_min}.
Bisection using \(\tau_D\) gives
\(\beta_{\min}(1)\approx0.091\) and
\(\beta_{\min}(2)\approx0.045\).
These predictions use only the centered additive
write-noise model; they exclude input/output noise,
quantization, and weight range variation.
Figure~\ref{fig:bcsstk21} compares these thresholds
with the measured iteration counts, where the applied preconditioner includes all noise sources, quantization, and weight range variation.

Figure~\ref{fig:bcsstk21} shows that small shifts
are insufficient for the tested analog configurations.
All \(10\) runs fail to converge within \(1000\)
iterations at the tested \(\beta\leq0.03\)
for \(D=1\) and at \(\beta=0.01\) for \(D=2\).
Successful solves first appear at \(\beta=0.045\)
and \(0.02\), respectively, below the predicted
thresholds.
Near the thresholds, the iteration counts fall
to \(452.7\pm7.6\) at \(\beta=0.09\) for \(D=1\)
and \(415.8\pm8.3\) at \(\beta=0.045\) for \(D=2\).

Further regularization initially improves performance.
The lowest observed analog counts are \(233.3\pm7.0\)
at \(\beta=0.35\) for \(D=1\) and \(129.6\pm7.7\)
at \(\beta=0.2\) for \(D=2\), compared with
\(100.1\pm0.9\) for the digital reference.
Larger shifts then increase the iteration counts.
This behavior illustrates why improving relative
write-error tolerance alone does not determine solver
performance: the shift also changes the inverse
approximation, as shown in
\eqref{eq:regularization_error}.

Copy averaging moves both the predicted threshold
and the best tested shift toward smaller values
of \(\beta\).
However, successful solves occur below
\(\beta_{\min}(D)\), while the fewest iterations
occur above it.
The threshold therefore provides a reference for
controlling local write error, not a convergence
boundary or an iteration-optimal shift for the
complete preconditioner.

The criterion in \eqref{eq:global_survival_criterion} is also consistent with an earlier observation that
regularization did not improve analog block Jacobi for the 3D
Poisson problems of~\cite{shah2026hybrid}: for those blocks, $\varphi(0) < \tau_1$,
so $\beta_{\min} = 0$ and no shift was needed to keep each block response above the write-error level.

\section{Conclusion}
\label{sec:conclusion}

We developed an analog--digital preconditioner combining
analog block Jacobi with a digital low-rank residual
correction learned from forward probes of the programmed
copies. Under explicit noise models, the analysis separates
programmed-mean error from read fluctuations, quantifies
probing accuracy, and provides guidance for averaging and
local regularization. In simulator tests with significant
programming errors, the learned spaces reduce FGMRES
iterations relative to digitally constructed spaces of
the same rank using the same analog base. The regularization
threshold controls local relative write error but does not
predict the shift giving the fewest iterations.

These results support constructing the correction from
the programmed action rather than the digital target alone.
End-to-end hardware evaluation remains necessary to determine
whether the iteration savings justify the setup, digital
correction, and data-movement costs. Learning $Z$ improves convergence at fixed $D$ and $S$, whereas averaging
improves it at the cost of additional arrays or reads; choosing $(k,D,S)$
jointly requires an architecture-specific energy and latency model, which we
leave for future work. Finally, a natural next step is a spectral analysis of Krylov methods with random,
zero-mean perturbations of a fixed preconditioner, which would relate the
one-step quantity $\gamma_Z^2(W) + \widehat{\sigma}_Z^2(W)/(DS)$ to FGMRES
convergence.

\section*{Acknowledgments}
We thank Tayfun Gokmen for productive conversations on this work. The authors gratefully acknowledge the support of the U.S. Department of Energy, Office of Science under grant DE-AC52-07NA27344. Claude Opus 5.5 and GPT-6 Astra were used to assist with literature reviews, producing code for figures, and the final preparation of the manuscript. The authors assume responsibility for all content.

\bibliographystyle{siamplain}
\bibliography{reference}
\end{document}